 \RequirePackage{fix-cm}

 \documentclass[smallextended]{svjour3} 
 \smartqed  
\usepackage{amssymb,amsmath,amsfonts,latexsym} 
\usepackage{hyperref}
\usepackage{bm}
\usepackage{algorithm}
\usepackage{algpseudocode}
\usepackage{tabularx}
\usepackage{booktabs}
 \def\proofend{\hfill$\Box$} 

\usepackage{mathtools}
\usepackage{multirow,array}
\usepackage{graphicx}
\usepackage{subfigure}
\usepackage{epstopdf}
\usepackage{float} 
\usepackage{color, xcolor}
\allowdisplaybreaks

\journalname{To appear in Fract. Calc. Appl. Anal.}

\begin{document}


\title{A fast fractional block-centered finite difference method for two-sided space-fractional diffusion equations on general nonuniform grids}

\titlerunning{A fast fractional block-centered finite difference method for \dots}
\author{
        Meijie Kong$^1$ 
\and
        Hongfei Fu$^2$ 
 }
\authorrunning{M. Kong \and  H. Fu} 
\institute{Meijie Kong$^{1}$
\at
School of Mathematical Sciences, Ocean University of China, Qingdao, Shandong 266100, China \\
\email{kmj@stu.ouc.edu.cn} 
 \and
Hongfei Fu$^{2,*}$
\at
School of Mathematical Sciences \&  Laboratory of Marine Mathematics, Ocean University of China, Qingdao, Shandong 266100, China. \\
\email{fhf@ouc.edu.cn} $^*$ corresponding author 
}

\date{Received: 11 December 2023 / Revised: 18 September 2024 / Accepted: ......}

\maketitle

\begin{abstract}
In this paper, a two-sided variable-coefficient space-fractional diffusion equation with fractional Neumann boundary condition is considered. To conquer the weak singularity caused by nonlocal space-fractional differential operators, a fractional block-centered finite difference (BCFD) method on general nonuniform grids is proposed. However, this discretization still results in an unstructured dense coefficient matrix with huge memory requirement and computational complexity.  To address this issue, a fast version fractional BCFD algorithm by employing the well-known sum-of-exponentials (SOE) approximation technique is also proposed. Based upon the Krylov subspace iterative methods, fast matrix-vector multiplications of the resulting coefficient matrices with any vector are developed, in which they can be implemented in only $\mathcal{O}(MN_{exp})$ operations per iteration without losing any accuracy compared to the direct solvers, where $N_{exp}\ll M$ is the number of exponentials in the SOE approximation. Moreover, the coefficient matrices do not necessarily need to be generated explicitly, while they can be stored  in $\mathcal{O}(MN_{exp})$ memory by only storing some coefficient vectors. Numerical experiments are provided to demonstrate the efficiency and accuracy of the method.



\end{abstract} 

	\section{Introduction}\label{sec:Intro}
	
	\setcounter{section}{1} \setcounter{equation}{0} 
	
	Fractional partial differential equations provide a very adequate and competitive tool to model challenging phenomena involving anomalous diffusion or long-range memory and spatial interactions. For example, the space-fractional diffusion equations (SFDEs) can be used to describe the anomalous diffusion occurred in many transport processes \cite{BWM00,SLHS15,MK00,SZCR14}. 
	Extensive research has been conducted in the development of numerical methods for various SFDEs \cite{DJWWZ21,ER06,FZS22,HPLC19,LZWF22,LFWC19,MS16,PSXX21,SCZC19,TZD15,XLS22,ZLLB14}. However, due to the nonlocal nature of the space-fractional differential operators, numerical discretizations tend to generate dense and full stiffness matrices. Traditionally, these methods are solved via the direct solvers such as Gaussian elimination (GE) method, which requires a computational complexity of order $\mathcal{O}(M^3)$ per time level and memory of order $\mathcal{O}(M^2)$, where $M$ is the total number of spatial unknowns in the numerical discretization. Consequently, numerical simulations of SFDEs would lead to significantly increased memory requirement and computational complexity as $M$ increases. However, in the case of uniform spatial partitions, Toeplitz-like structure of the resulting stiffness matrices was  discovered in \cite{WWS10} for space-fractional diffusion model, and thus based upon the special matrix structures, fast Krylov subspace iterative solvers are developed for different numerical methods of various SFDEs, in which both memory requirement and computational complexity have been largely reduced \cite{FLW19,JW15fb,JW16,LZFW21,PKNS14,WB12,WD13,WW11,ZZFL21}. 
	
	On the other hand, it is shown that the solutions to SFDEs usually exhibit power function singularities near the boundary even under assumptions that the coefficients and right-hand side source term are sufficiently smooth \cite{EHR18,JLPR15,WY17}. This motivates the usage of nonuniform grids to better capture the singular behavior of the solutions. Simmons \cite{SYM17} derived a first-order finite volume method on nonuniform grids for two-sided fractional diffusion equations. However, since the dense coefficient matrices do not have Toeplitz-like structure as the case of  uniform grids, the previous developed fast algorithms are no longer be applicable. It makes sense to construct novel fast algorithms discretized on arbitrary nonuniform grids. Recently, finite volume methods based on special composite structured grids that consist of a uniform spatial partition far away the boundary and a locally refinement near the boundary were studied to address this issue \cite{DJWWZ21,JW15,JW19EAJAM,JW19}, where Toeplitz-like structures of the coefficient matrices can be still found, and the Toeplitz structure of the diagonal blocks of the resulting block coefficient matrix was employed for efficient evaluate matrix-vector multiplications and the off-diagonal blocks were properly approximated by low-rank decompositions. 
	However, there are still lack of fast numerical methods for SFDEs on general nonuniform grids. Recently, Jiang et al. \cite{JZZZ17} developed a novel sum-of-exponentials (SOE) technique for fast evaluation of the Caputo time-fractional derivative and applied to time-fractional diffusion equations. Then, by adopting the SOE technique on graded spatial grids, Fang et al. \cite{FZS22}  proposed a  fast first-order finite volume method for the one-dimensional SFDEs  with homogeneous Dirichlet boundary condition. However, to our best knowledge, currently there seems  rarely papers on construction of highly efficient numerical methods on general nonuniform spatial grids, and it is still a major challenge for the modeling of SFDEs.  
	
	Block-centered finite difference (BCFD) method, sometimes called cell-centered finite difference method, can be thought of as the lowest-order Raviart-Thomas mixed element method, by employing a proper numerical quadrature formula \cite{RT77}. One of the most important merits of the BCFD method is that it can simultaneously approximate the primal variable and its flux to a same order of accuracy on nonuniform grids, without any accuracy lost compared to the standard finite difference method. Besides, the BCFD method can  very easily deal with model problems with Neumann or periodic boundary conditions. In recent years, BCFD method has been widely used to simulate integer-order PDEs and time-fractional PDEs \cite{LL16,RP12,XXF22}. However, due to the complexity of the nonlocal space-fractional differential operators, we do not see any report on BCFD method for space-fractional PDEs. Therefore, we aim to present a fractional type BCFD method for model \eqref{model-2d:e1} with fractional Neumann boundary condition, in which nonuniform spatial grids are utilized to capture the boundary singularity of the solution and thus to improve the computational accuracy. Moreover, another main goal of this paper is to present a fast  fractional BCFD algorithm to improve the computational efficiency of  modeling the SFDEs, in which a fast Krylov subspace iterative solver  based upon efficient matrix-vector multiplications is developed. 
	
	The present paper focuses on efficient and accurate numerical approximation of the following two-sided variable-coefficient SFDE with anomalous diffusion orders $1<\alpha,\beta<2$ \cite{SBMW01}:
	\begin{equation}\label{model-2d:e1}
		\begin{aligned}
			\partial_t u - \partial_x \, \mathcal{D}_{\gamma}^{x,\alpha}u(x,y,t)&- \partial_ y\, \mathcal{D}_\theta^{y,\beta}u(x,y,t)=f(x,y,t),\\
			& \qquad \quad (x,y)\in\Omega:=(a,b)\times(c,d),\quad 0<t\leq T,
		\end{aligned}
	\end{equation}
	where $\mathcal{D}_{\gamma}^{x,\alpha}$ represent the weighted fractional-order differential operators along $x$ direction:
	\begin{equation*}
		\begin{aligned}
			\mathcal{D}_{\gamma}^{x,\alpha}&:=\gamma K^{x,L} \partial_x \, {}_a\mathcal{J}_x^{2-\alpha} +(1-\gamma)K^{x,R} \partial_x \, {}_x\mathcal{J}_b^{2-\alpha}, 
		\end{aligned}
	\end{equation*}
	with weighted parameters $0\leq\gamma\leq 1$ is a parameter describing the relative probability of particle traveling ahead or behind the mean velocity and positive diffusion coefficients $K^{x,L}=K^{x,L}(x,y,t)$, $K^{x,R}=K^{x,R}(x,y,t)$. Physically, the fractional derivative $\mathcal{D}_{\gamma}^{x,\alpha}$ can be interpreted as a nonlocal Fickian law \cite{SBMW01}. Moreover, ${}_a\mathcal{J}_x^{2-\alpha}$ and ${}_x\mathcal{J}_b^{2-\alpha}$ respectively denote the left/right-sided Riemann-Liouville fractional integrals, defined by \cite{P99} 
	\begin{equation*}
		\begin{aligned}
			{ }_a\mathcal{J}_x^\alpha g(x,y)&:=\frac{1}{\Gamma(\alpha)}\int_a^x(x-s)^{\alpha-1} g(s,y)ds,\\
			{ }_x\mathcal{J}_b^\alpha g(x,y)&:=\frac{1}{\Gamma(\alpha)}\int_x^b(s-x)^{\alpha-1} g(s,y)ds.
	    \end{aligned}
	 \end{equation*}
Similarly, the operator $\mathcal{D}_\theta^{y,\beta}$ along $y$ direction can be also defined.

	We assume that problem \eqref{model-2d:e1} is subjected to the following initial condition
	\begin{equation}\label{model-2d:e2}
		u(x,y,0)=u^o(x,y),\quad\text{for}\quad (x,y)\in\bar{\Omega},
	\end{equation}
	and fractional Neumann boundary conditions \cite{BKMSS15,JW15fb,WWC18}
	\begin{equation}\label{model-2d:e3}
		\begin{aligned}
			&\mathcal{D}_{\gamma}^{x,\alpha} u(x,y,t)\mid_{x=a} =\phi^x(y,t), \quad\mathcal{D}_{\gamma}^{x,\alpha} u(x,y,t)\mid_{x=b}=\varphi^x(y,t),\\
			&\mathcal{D}_\theta^{y,\beta} u(x,y,t)\mid_{y=c}=\phi^y(x,t),\quad\mathcal{D}_\theta^{y,\beta} u(x,y,t)\mid_{y=d}=\varphi^y(x,t).
		\end{aligned}
	\end{equation}
	
In this paper,  by introducing an auxiliary fractional flux variable, we first develop a fractional CN-BCFD method  on general nonuniform grids for SFDEs \eqref{model-2d:e1}--\eqref{model-2d:e3} in one space dimension, in which the Crank-Nicolson (CN) temporal discretization combined with the fractional BCFD spatial discretization are employed. 
	 Then, we develop fast approximation techniques to the left/right-sided Riemann-Liouville fractional integrals, using efficient SOE approximations to the kernels $x^{1-\alpha}$, $y^{1-\beta}$ and piecewise linear interpolations of the primal variable $u$. Based upon these approximations, fast Krylov subspace iterative solvers for the  fractional BCFD method is then proposed with fast matrix-vector multiplications of the resulting coefficient matrices with any vector. It is shown that the solver requires only $\mathcal{O}(MN_{exp})$ operations per iteration with efficient matrix storage mechanism, where $N_{exp}\ll M$ is the number of exponentials in the SOE approximation. Finally, ample numerical experiments are provided to demonstrate the efficiency and  accuracy of the method.
	As far as we know, this seems to be the first time that a fast fractional BCFD method is developed for the SFDEs, where fractional Neumann boundary condition is considered and  general nonuniform spatial grids are adopted.	
	
	The rest of the paper is organized as follows. In Section \ref{sec:dbcfd}, we present the fractional CN-BCFD method on nonuniform spatial grids for the SFDEs. In Section \ref{sec:fbcfd}, a fast version fractional CN-BCFD method is proposed to further improve the computational efficiency. Then, we give some numerical experiments to investigate the accuracy and  performance of the fast method in Section \ref{sec:num}. Some concluding remarks are given in the last section.
	
	\section{A fractional CN-BCFD method on general nonuniform grids}\label{sec:dbcfd}
	\setcounter{section}{2} \setcounter{equation}{0} 
	
	For simplicity of presentation, in this section we pay our attention to the  one-dimensional version of  \eqref{model-2d:e1}:
	\begin{equation}\label{model-1d}
		\left\lbrace 
		\begin{aligned}{}
			& \partial_t u-\partial_x \mathcal{D}_{\gamma}^{x,\alpha}u(x,t)=f(x,t),\quad  x\in(a,b),\ 0<t\leq T,\\
			&  \mathcal{D}_{\gamma}^{x,\alpha}u(x,t)\big|_{x=a} =\phi(t),\,\mathcal{D}_{\gamma}^{x,\alpha}u(x,t)\big|_{x=b}=\varphi(t), \quad 0\leq t\leq T,\\
			& u(x,0)=u^o(x),\quad x\in[a,b],
		\end{aligned}
		\right. 
	\end{equation}
with $\mathcal{D}^{x,\alpha}_{\gamma}=\gamma K^L\partial_x\, {}_a\mathcal{J}_x^{2-\alpha}+(1-\gamma)K^R\partial_x\,{}_x\mathcal{J}_b^{2-\alpha}$. While, the extension to the two-dimensional case is  straightforward but more complicated.
	
	In the following, we introduce some notations.First, for positive integer $N$, we define a uniform temporal partition of $[0,T]$ by $t_n=n\tau$, $n=0,1,\cdots,N$, with  temporal stepsize $\tau:=T/N$. For temporal grid function $\left\{\Phi^n\right\}_{n=0}^N$, define
	\[
	\Phi^{n-1/2}:=\frac{\Phi^n+\Phi^{n-1}}{2},\quad\delta_t\Phi^n:=\frac{\Phi^n-\Phi^{n-1}}{\tau}.
	\]
Next, for positive integer $M$, define a set of nonuniform  spatial grids by
	\[
	a=x_{1/2}<x_{3/2}<\cdots<x_{M-1/2}<x_{M+1/2}=b,
	\]
	with grid size $h_i=x_{i+1/2}-x_{i-1/2}$ for $i=1,2,\dots,M$ and $h:=\max_{1\leq i\leq M}h_i$. 
Besides, set $x_i:=(x_{i+1/2}+x_{i-1/2})/2$ for $i=1,2,\ldots,M$ as another set of nonuniform staggered spatial grids. Denote $h_{1/2}:=x_1-a=h_1/2$, $h_{i+1/2}:=x_{i+1}-x_i= (h_{i+1}+h_i)/2$, $1\leq i\leq M-1$, $h_{M+1/2}:=b-x_{M}= h_{M}/2$. Furthermore, for spatial grid functions $\left\{\psi_i\right\}_{i=1}^{M}$ and $\left\{\psi_{i+1/2}\right\}_{i=0}^{M}$, define
	\[
	d_x\psi_{i+1/2}:=\frac{\psi_{i+1}-\psi_i}{h_{i+1/2}},\quad D_x\psi_i:=\frac{\psi_{i+1/2}-\psi_{i-1/2}}{h_i}.
	\]
	
	To propose a second-order fractional BCFD method for model \eqref{model-1d}, we introduce the piecewise linear interpolation function $\mathcal{L}_h[v](x)$ for $v(x)$ as follows:
	\begin{equation}\label{equ:interpolation}
		\mathcal{L}_h[v](x)= \left\{
		\begin{aligned}
			&\frac{x_1-x}{h_{1/2}}\bar{v}(a)+\frac{x-a}{h_{1/2}}v_1,  \quad x\in [a,x_1],\\ 
			&\frac{x_{i+1}-x}{h_{i+1/2}}v_i+\frac{x-x_i}{h_{i+1/2}}v_{i+1},    ~ x\in[x_i,x_{i+1}],\ 1\leq i\leq M-1,\\ 
			&\frac{b-x}{h_{M+1/2}}v_{M}+\frac{x-x_{M}}{h_{M+1/2}}\bar{v}(b),  \quad x\in [x_{M},b],
		\end{aligned}\right.
	\end{equation}
	where $v_i:=v(x_i)$, and $\bar{v}(a)$ and $\bar{v}(b)$ are defined by two-point extrapolations respectively as
	\[
	\bar{v}(a):=\frac{\left(2h_1+h_2\right)v_1-h_1v_2}{h_1+h_2},\quad\bar{v}(b):=\frac{\left(2h_{M}+h_{M-1}\right)v_{M}-h_{M}v_{M-1}}{h_{M-1}+h_{M}}.
	\]
	By Taylor's expansion, it is easily to show that for smooth $v$ it holds
	\[
	\left|(\bar{v}-v)(a)\right|=O(h^2),\quad\left|(\bar{v}-v)(b)\right|=O(h^2).
	\]
	
	In this section, we aim to develop a direct fractional BCFD method on general nonuniform spatial grids for the one-dimensional SFDE \eqref{model-1d}. To this aim, we introduce an auxiliary fractional flux variable
	\begin{equation}\label{model:flux}
		\begin{aligned}
			p(x,t)=\gamma K^{L}\,\partial_x \, g^{L}(x,t)+(1-\gamma)K^{R}\,\partial_x \, g^{R}(x,t),
		\end{aligned}
	\end{equation}
	with
	\begin{equation}\label{model:flux:e2}
		\begin{aligned}
			g^{L}(x,t)={}_a\mathcal{J}_x^{2-\alpha}u(x,t),\quad
			g^{R}(x,t)={}_x\mathcal{J}_b^{2-\alpha}u(x,t).
		\end{aligned}
	\end{equation}
	Then, the original two-sided variable-coefficient SFDE  model \eqref{model-1d} is equivalent to
	\begin{equation}\label{equ:rewrite-model-1d}
		\left\lbrace 
		\begin{aligned}{}
			& \partial_t u-\partial_x  p(x,t)=f(x,t),\quad x\in(a,b),\ 0<t\leq T, \\
			& p(a,t)=\phi(t),\quad p(b,t)=\varphi(t), \quad 0<t\leq T,\\
			& u(x,0)=u^o(x),\quad x\in[a,b],
		\end{aligned}
		\right. 
	\end{equation}
	where $p(x,t)$ is given by \eqref{model:flux}--\eqref{model:flux:e2}.
	
	The BCFD method can be thought of as a special mixed element method, in which on each element $[x_{i-1/2},x_{i+1/2}]$, the flux variable $p$ is approximated at the endpoint of each element, i.e., $\{x_{i+1/2}\}$, while the primal variable $u$ is approximated at the midpoint of each element, i.e., $\{ x_{i}\}$.	Denote the finite difference approximations $u_i^n\approx u(x_i,t_n)$ for $1\leq i\leq M$, $1\leq n\leq N$ and $p_{i+1/2}^{n}\approx p(x_{i+1/2},t_n)$ for $0\leq i\leq M$, $1\leq n\leq N$. Then,  a second-order Crank-Nicolson semi-discretization scheme for \eqref{equ:rewrite-model-1d} reads as
	\begin{equation}\label{equ:time-semi-discrete}
		\left\{
		\begin{aligned}{}
			& \delta_tu^n(x)- \frac{1}{2} \partial_x  p(x,t_{n}) = \frac{1}{2}\partial_x  p(x,t_{n-1})+		f(x,t_{n-1/2}), \  x\in(a,b), \\
			&p(a,t_n)=\phi(t_n),\quad  p(b,t_n)=\varphi(t_n), \\
			& u(x,0)=u^o(x),\quad x\in[a,b],
		\end{aligned}
		\right. 
	\end{equation}
	for $1 \le n\leq N$, where $p(x,t_n)$ is defined by \eqref{model:flux}--\eqref{model:flux:e2} at time $t=t_n$.
	
We next consider the spatial discretization of \eqref{equ:time-semi-discrete} using the BCFD method on staggered grids  $\{x_{i}\}$ and $\{x_{i+1/2}\}$. The crucial step is the approximations to the fractional integrals $g^L$ and $g^R$ in \eqref{model:flux:e2} at each grid point $x=x_i$. To this aim, we first  construct an approximation to the left-sided Riemann-Liouville fractional integral $g^L$. By splitting the integral  into two parts, of which the first part is an integral on half a grid interval, and then approximating the unknown function $u(x,t)$ by its  linear interpolation \eqref{equ:interpolation},  we get 
	\begin{equation}\label{equ:left-inte-appr}
		\begin{aligned}
		&\quad 	g^{L}(x_i,t_n) \\
			& \approx\int_a^{x_1}\frac{(x_i-\xi)^{1-\alpha}}{\Gamma(2-\alpha)}\mathcal{L}_h[u^n](\xi)d\xi+ \sum_{j=1}^{i-1}\int_{x_j}^{x_{j+1}}\frac{(x_i-\xi)^{1-\alpha}}{\Gamma(2-\alpha)}\mathcal{L}_h[u^n](\xi) d\xi\\
			&=\int_a^{x_1}\frac{(x_i-\xi)^{1-\alpha}}{\Gamma(2-\alpha)} \left[\frac{x_1-\xi}{h_{1/2}}\bar{u}^n(a)+\frac{\xi-a}{h_{1/2}}u_1^n\right]d\xi\\
			& \qquad +\sum_{j=1}^{i-1}\int_{x_j}^{x_{j+1}}\frac{(x_i-\xi)^{1-\alpha}}{\Gamma(2-\alpha)} \left[\frac{x_{j+1}-\xi}{h_{j+1/2}}u_j^n+\frac{\xi-x_j}{h_{j+1/2}}u_{j+1}^n\right]d\xi\\
			&=\sum_{j=1}^{i}q_{i,j}^Lu_j^n=:g_i^{L,n},
		\end{aligned}
	\end{equation}
for $2 \le i \le M$, and, in particular, for $i=1$, \eqref{equ:left-inte-appr} reduces to
	\begin{equation}\label{equ:left-inte-appr1}
		g^{L}(x_1,t_n)\approx\int_a^{x_1}\frac{(x_1-\xi)^{1-\alpha}}{\Gamma(2-\alpha)}\mathcal{L}_h[u^n](\xi)d\xi
		=q_{1,1}^Lu_1^n+q_{1,2}^Lu_2^n=: g_{1}^{L,n}.
	\end{equation}
	A simple calculation shows that the coefficients in \eqref{equ:left-inte-appr}--\eqref{equ:left-inte-appr1} can be expressed as
	\begin{equation}\label{equ:coefficient-ql}\left\{
		\begin{aligned}
			q^L_{i,1}&=\frac{2h_1+h_2}{h_1+h_2}\frac{(x_i-a)^{2-\alpha}}{\Gamma(3-\alpha)} +\frac{h_1}{h_1+h_2}\omega_{i,1}^L+\omega_{i,2}^L, \\
			q^L_{i,2}& =-\frac{h_1}{h_1+h_2}\frac{(x_i-a)^{2-\alpha}}{\Gamma(3-\alpha)} -\frac{h_1}{h_1+h_2}\omega_{i,1}^L-\omega_{i,2}^L+\omega_{i,3}^L,\\
			q^L_{i,j}&=\omega_{i,j+1}^L-\omega_{i,j}^L,\quad 3\leq j\leq i,
		\end{aligned}\right.
	\end{equation}
	with
	\begin{equation}\label{equ:coefficient-wl}\left\{
		\begin{aligned}
			&\omega_{i,1}^L=\frac{(x_i-x_1)^{3-\alpha}-(x_i-a)^{3-\alpha}}{\Gamma(4-\alpha)h_{1/2}},\\
			&\omega_{i,j}^L=\frac{(x_i-x_j)^{3-\alpha}-(x_i-x_{j-1})^{3-\alpha}}{\Gamma(4-\alpha)h_{j-1/2}}, \quad  2\leq j\leq i,\\ 
			& \omega_{i,j}^L=0,\quad j>i,
		\end{aligned}\right.
	\end{equation}
	for $1\leq i\leq M$.
	
	Next, we pay attention to the  approximation to the right-sided  Riemann-Liouville fractional integral $g^R$.	Similarly,  for $1\leq i\leq M-1$, it can be approximated using the linear interpolation \eqref{equ:interpolation} by
\begin{equation}\label{equ:right-inte-appr}
		\begin{aligned}
		&\quad 	g^{R}(x_i,t_n)\\
			& \approx\sum_{j=i}^{M-1}\int_{x_j}^{x_{j+1}}\frac{(\xi-x_i)^{1-\alpha}}{\Gamma(2-\alpha)}\mathcal{L}_h[u^n](\xi)d\xi+\int_{x_M}^{b}\frac{(\xi-x_i)^{1-\alpha}}{\Gamma(2-\alpha)}\mathcal{L}_h[u^n](\xi)d\xi\\
			& =\sum_{j=i}^{M-1}\int_{x_j}^{x_{j+1}}\frac{(\xi-x_i)^{1-\alpha}}{\Gamma(2-\alpha)} \left[\frac{x_{j+1}-\xi}{h_{j+1/2}}u_j^n+\frac{\xi-x_j}{h_{j+1/2}}u_{j+1}^n\right]d\xi\\
			&~ \qquad +\int_{x_M}^{b}\frac{(\xi-x_i)^{1-\alpha}}{\Gamma(2-\alpha)} \left[\frac{b-\xi}{h_{M+1/2}}u_M^n+\frac{\xi-x_M}{h_{M+1/2}}\bar{u}^n(b)\right]d\xi\\
			&  =\sum_{j=i}^Mq^R_{i,j}u_j^n=:g_i^{R,n},
		\end{aligned}
	\end{equation}
	and, in particular, for $i=M$, \eqref{equ:right-inte-appr} reduces to
	\begin{equation}\label{equ:right-inte-apprM}
		g^{R}(x_M,t_n)\approx\int_{x_M}^{b}\frac{(\xi-x_M)^{1-\alpha}}{\Gamma(2-\alpha)}\mathcal{L}_h[u^n](\xi) d\xi=q^R_{M-1,M}u_{M-1}^n+q^R_{M,M}u_M^n=:g_M^{R,n}.
	\end{equation}
	Similarly, the coefficients in \eqref{equ:right-inte-appr}--\eqref{equ:right-inte-apprM} are expressed as
	\begin{equation}\label{equ:coefficient-qr}\left\{
		\begin{aligned}
			&q^R_{i,j}=\omega_{i,j}^R-\omega_{i,j-1}^R,\quad i\leq j\leq M-2,\\
			&q^R_{i,M-1}= \frac{-h_M}{h_{M-1}+h_M}\frac{(b-x_i)^{2-\alpha}}{\Gamma(3-\alpha)}
			-\omega_{i,M-2}^R+\omega_{i,M-1}^R
			+\frac{h_M}{h_{M-1}+h_M}\omega_{i,M}^R,\\ 
			& q_{i,M}^R=\frac{2h_M+h_{M-1}}{h_M+h_{M-1}}\frac{(b-x_i)^{2-\alpha}}{\Gamma(3-\alpha)}
			-\omega_{i,M-1}^R-\frac{h_M}{h_{M-1}+h_M}\omega_{i,M}^R,
		\end{aligned}\right.
	\end{equation}	
	with
	\begin{equation}\label{equ:coefficient-wr}\left\{
		\begin{aligned}
			&\omega_{i,j}^R=\frac{(x_{j+1}-x_i)^{3-\alpha}-(x_{j}-x_i)^{3-\alpha}}{\Gamma(4-\alpha)h_{j+1/2}},\quad    i\leq j\leq M-1,\\
			&\omega_{i,M}^R=\frac{(b-x_i)^{3-\alpha}-(x_{M}-x_i)^{3-\alpha}}{\Gamma(4-\alpha)h_{M+1/2}},\\
			&\omega_{i,j}^R=0,\quad i>j,
		\end{aligned}\right.
	\end{equation}
	for $1\leq i\leq M$.
	
	Now,  by combining \eqref{equ:left-inte-appr}--\eqref{equ:left-inte-appr1} and 	\eqref{equ:right-inte-appr}--\eqref{equ:right-inte-apprM} with \eqref{equ:time-semi-discrete}, a fully discrete fractional CN-BCFD scheme  is proposed as follows:
	\begin{equation}\label{equ:sche-cn-bcfd}
		\left\{
		\begin{aligned}{}
			&\delta_tu_i^n- \frac{1}{2} D_x p_i^{n}	=\frac{1}{2} D_x p_i^{n-1}	+f_i^{n-1/2}, \quad 1\leq i\leq M,\\
			& p_{i+1/2}^{n}=\gamma K^{L,n}_{i+1/2}\ d_xg^{L,n}_{i+1/2}+(1-\gamma) K^{R,n}_{i+1/2}\ d_x g^{R,n}_{i+1/2}, ~ 1\leq i\leq M-1,\\
			& p^{n}_{1/2}=\phi^{n},\quad	 p^{n}_{M+1/2}=\varphi^{n}, \\
			& u_i^0=u^o(x_i),\quad 1\leq i\leq M,
		\end{aligned}
		\right.
	\end{equation}
	where $K^{L,n}_{i+1/2}:=K^{L}(x_{i+1/2},t_n)$, $K^{R,n}_{i+1/2}:=K^{R}(x_{i+1/2},t_n)$, $\phi^{n}:=\phi(t_n)$ and $\varphi^{n}:=\varphi(t_n)$.
	Moreover, let
	\begin{equation*} 
		\begin{aligned}
			\bm{u}^n &=[u_1^n,u_2^n,\dots,u_M^n]^\top,\\
			\bm{F}^{n-1/2}&=\Big[f_1^{n-1/2}-\frac{\phi^{n}+ \phi^{n-1}}{2h_1},f_2^{n-1/2},\dots,f_{M-1}^{n-1/2},
			f_M^{n-1/2}+\frac{\varphi^{n}+ \varphi^{n-1}}{2h_M}\Big]^\top,
		\end{aligned}
	\end{equation*}
	where $\top$ refers to the transpose of the vector.  Then, by canceling the flux variable, we can present  the fractional CN-BCFD scheme \eqref{equ:sche-cn-bcfd} in a more compact  matrix form  with respect to $\bm{u}^n$:
	\begin{equation}\label{equ:matrix form}
		\left(\bm{I}_M-\frac{\tau}{2}\bm{A}^n\right)\bm{u}^n=\left(\bm{I}_M+\frac{\tau}{2}\bm{A}^{n-1}\right)\bm{u}^{n-1}+\tau\bm{F}^{n-1/2},
	\end{equation}
	where $\bm{I}_M$ represents the identity matrix of order $M$, 	and $\bm{A}^n$ is a stiffness matrix of order $M$ with entries corresponding to $D_x p_i^{n}$ such that
	\begin{equation}\label{coeff-a}
		\begin{aligned}
		\bm{A}^n& =\gamma\left[ \textrm{diag}\big(\bm{D}_+^{L,n}\big)\bm{A}_+^{L}-\textrm{diag} \big(\bm{D}_-^{L,n}\big) \bm{A}_-^L\right] \\
		&\qquad +(1-\gamma)\left[ \textrm{diag}\big( \bm{D}_+^{R,n}\big) \bm{A}_+^{R}-\textrm{diag} \big(\bm{D}_-^{R,n}\big) \bm{A}_-^R\right],
		\end{aligned}
	\end{equation}
	where,  for simplicity, hereafter we denote $\mathcal{K}^{L,n}_{i+1/2}:=K^{L,n}_{i+1/2}/{h_{i+1/2}}$ and $ \mathcal{K}^{R,n}_{i+1/2}:={K^{R,n}_{i+1/2}}/{h_{i+1/2}}$, and
	\begin{equation}\label{coeff:e1}
		\bm{D}_+^{L,n}=\Big[ \frac{\mathcal{K}^{L,n}_{3/2}}{h_1},\cdots,\frac{\mathcal{K}^{L,n}_{M-1/2}}{h_{M-1}},0\Big]^\top,
		\ 
		\bm{D}_-^{L,n}=\Big[ 0,\frac{\mathcal{K}^{L,n}_{3/2}}{h_2},\cdots,\frac{\mathcal{K}^{L,n}_{M-1/2}}{h_{M}}\Big]^\top,
	\end{equation}
	\begin{equation}\label{coeff:e2}
		\bm{D}_+^{R,n}=\Big[ \frac{\mathcal{K}^{R,n}_{3/2}}{h_1},\cdots,\frac{\mathcal{K}^{R,n}_{M-1/2}}{h_{M-1}},0\Big]^\top,
		\
		\bm{D}_-^{R,n}=\Big[ 0,\frac{\mathcal{K}^{R,n}_{3/2}}{h_2},\cdots,\frac{\mathcal{K}^{R,n}_{M-1/2}}{h_{M}}\Big]^\top,
	\end{equation}
	and
	\begin{equation*}
		\begin{aligned}
			\bm{A}_+^L=\begin{bmatrix}
				q_{2,1}^L-q_{1,1}^L   & q_{2,2}^L-q_{1,2}^L   & \cdots & 0  & 0         \\
				q_{3,1}^L-q_{2,1}^L   & q_{3,2}^L-q_{2,2}^L   & \ddots & \ddots   & 0         \\
				q_{4,1}^L-q_{3,1}^L   & q_{4,2}^L-q_{3,2}^L   & \ddots & \ddots & \vdots    \\
				\vdots                & \vdots                &        & q_{M-1,M-1}^L   & 0         \\ 
				q_{M,1}^L-q_{M-1,1}^L & q_{M,2}^L-q_{M-1,2}^L & \cdots & q_{M,M-1}^L-q_{M-1,M-1}^L & q_{M,M}^L \\ 
				0                     & 0                     & \cdots & 0  & 0         \\
			\end{bmatrix},
		\end{aligned}
	\end{equation*}
	\begin{equation*}
		\begin{aligned}
			\bm{A}_-^L=\begin{bmatrix}
				0                     & 0                     & \cdots & 0                         & 0      \\
				q_{2,1}^L-q_{1,1}^L   & q_{2,2}^L-q_{1,2}^L   & \cdots & 0                         & 0      \\
				q_{3,1}^L-q_{2,1}^L   & q_{3,2}^L-q_{2,2}^L   & \ddots & \ddots                    & \vdots \\
				\vdots                & \ddots                & \ddots & \ddots                    & 0      \\ 
				q_{M,1}^L-q_{M-1,1}^L & q_{M,2}^L-q_{M-1,2}^L & \cdots & q_{M,M-1}^L-q_{M-1,M-1}^L & q_{M,M}^L\\ 
			\end{bmatrix},
		\end{aligned}
	\end{equation*}
	\begin{equation*}
		\begin{aligned}
			\bm{A}_+^R=
			\begin{bmatrix}
				-q^R_{1,1} & q^R_{2,2}-q^R_{1,2} & \cdots & q^R_{2,M-1}-q^R_{1,M-1} & q^R_{2,M}-q^R_{1,M}         \\
				0         & -q^R_{2,2}         & \cdots & q^R_{3,M-1}-q^R_{2,M-1} & q^R_{3,M}-q^R_{2,M}    \\
				\vdots    & \ddots            & \ddots & \vdots                  & \vdots               \\
				0         & 0                 & \cdots & q_{M,M-1}^R-q^R_{M-1,M-1}          & q^R_{M,M}-q^R_{M-1,M} \\
				0         & 0                 & \cdots & 0                       & 0                    \\
			\end{bmatrix},
		\end{aligned}
	\end{equation*}
	\begin{equation*}
		\begin{aligned}
			\bm{A}_-^R=
			\begin{bmatrix} 
				0         & 0                 & \cdots & 0                       & 0                    \\
				-q^R_{1,1} & q^R_{2,2}-q^R_{1,2} & \cdots & q^R_{2,M-1}-q^R_{1,M-1} & q^R_{2,M}-q^R_{1,M}    \\
				0         & -q^R_{2,2}         & \cdots & q^R_{3,M-1}-q^R_{2,M-1} & q^R_{3,M}-q^R_{2,M}    \\
				\vdots    & \ddots            & \ddots & \vdots                  & \vdots               \\
				0         & 0                 & \cdots & q_{M,M-1}^R-q^R_{M-1,M-1}          & q^R_{M,M}-q^R_{M-1,M} \\
			\end{bmatrix}.
		\end{aligned}
	\end{equation*}
	
	\begin{remark}
		It can be observed that the resulting stiffness matrix $\bm{A}^n$ is actually a full and dense matrix with a complicated structure, which requires $\mathcal{O}(M^2)$ memory for storage. Furthermore, at each time level, traditional Gaussian type direct solvers require $\mathcal{O}(M^3)$ computational complexity, and Krylov subspace iterative solvers also require $\mathcal{O}(M^2)$ computational complexity per iteration.	Therefore, the implementation of  \eqref{equ:matrix form} is indeed computationally expensive, especially for large-scale modeling and simulation. And thus, an efficient solution method for  \eqref{equ:matrix form} is of course highly demanded. We shall discuss this issue in next section.
	\end{remark}
	\begin{remark}
		In the case of uniform spatial partition, i.e., $h_i \equiv h$, the stiffness matrix $\bm{A}^n$ in \eqref{coeff-a} can be expressed as 
		\begin{equation}\label{coeff-b}
			\bm{A}^n=\gamma\bm{D}^{L,n}\bm{A}^L+(1-\gamma)\bm{D}^{R,n}\bm{A}^R,
		\end{equation}
		where
		\begin{equation*}
			\begin{split}
				\bm{D}^{L,n}=\frac{1}{h^2}\begin{bmatrix}
					0 & K^{L,n}_{3/2}&&&\\
					& -K^{L,n}_{3/2}&K_{5/2}^{L,n}&&\\
					&&\ddots &\ddots&\\
					&&&-K_{M-3/2}^{L,n}&K_{M-1/2}^{L,n}\\
					&&&&-K_{M-1/2}^{L,n}
				\end{bmatrix},
			\end{split}
		\end{equation*}
		\begin{equation*}
			\begin{split}
				\bm{D}^{R,n}=\frac{1}{h^2}\begin{bmatrix}
					K^{R,n}_{3/2}&&&\\
					-K^{R,n}_{3/2}& K^{R,n}_{5/2}&&&\\
					&\ddots&\ddots &&\\
					&&-K_{M-3/2}^{R,n}&K_{M-1/2}^{R,n}&\\
					&&&-K_{M-1/2}^{R,n}&0
				\end{bmatrix},
			\end{split}
		\end{equation*}
		and the matrices $\bm{A}^L$ and $\bm{A}^R$ can be expressed into the following $2\times 2$ blocks
		\begin{equation*}
			\begin{split}
				\bm{A}^{L}=\begin{bmatrix}
					\bm{A}_{2,2}^{L}& \bm{0}\\
					\bm{A}_{M-2,2}^{L}&\bm{A}_{M-2,M-2}^{L}
				\end{bmatrix},\quad\bm{A}^{R}=\begin{bmatrix}
					\bm{A}_{M-2,M-2}^{R}& \bm{A}_{2,M-2}^{R}\\
					\bm{0}&\bm{A}_{2,2}^{R}
				\end{bmatrix}.
			\end{split}
		\end{equation*}
		The submatrices $\bm{A}^L_{2,2}$ and $\bm{A}^R_{2,2}$ of order $2$ represent the stiffness matrices corresponding to the nodes $\left\lbrace x_1, x_2\right\rbrace $ and $\left\lbrace x_{M-1}, x_M\right\rbrace $, respectively. The submatrix $\bm{A}^L_{M-2,2}$ of order $(M-2)$-by-$2$ represents the coupling between the nodes $\left\lbrace x_3,\dots,x_M\right\rbrace $ and the nodes $\left\lbrace x_1, x_2\right\rbrace $, while $\bm{A}^R_{2,M-2}$ of order $2$-by-$(M-2)$ represents the coupling between the nodes $\left\lbrace x_1,\dots,x_{M-2}\right\rbrace $ and the nodes $\left\lbrace x_{M-1}, x_M\right\rbrace $. Finally, the submatrices $\bm{A}^L_{M-2,M-2}$ and $\bm{A}^R_{M-2,M-2}$ of order $M-2$ represent the stiffness matrices corresponding to the nodes $\left\lbrace x_3, \dots, x_M\right\rbrace $ and $\left\lbrace x_{1}, \dots, x_{M-2}\right\rbrace $, respectively, and both of them have special Toeplitz structures such that
		\begin{equation*}
			\begin{split}
				\bm{A}^L_{M-2,M-2}=\begin{bmatrix}
					q_1 &&&&\\
					q_2 & q_1 &&&\\
					q_3&q_2 & q_1 &&\\
					\vdots & \ddots & \ddots&\ddots&\\
					q_{M-2} & q_{M-3} &\cdots&q_2&q_1
				\end{bmatrix},
				\quad \bm{A}^R_{M-2,M-2}=\left(\bm{A}^L_{M-2,M-2}\right)^\top,
			\end{split}
		\end{equation*}
		where 
		\begin{equation*}
			q_i=\frac{h^{2-\alpha}}{\Gamma(4-\alpha)}\left\{\begin{aligned}
				&1,\quad i=1\\
				&2^{3-\alpha}-3,\quad i=2,\\
				&i^{3-\alpha}-3(i-1)^{3-\alpha}+3(i-2)^{3-\alpha}-(i-3)^{3-\alpha}, ~3\leq i\leq M-2.
			\end{aligned}
			\right.
		\end{equation*}
Thus,  the stiffness matrix $\bm{A}^n$ has Toeplitz-like structure, and therefore,  the resulting linear algebraic system \eqref{equ:matrix form} can be solved via the fast Fourier transform (FFT) approach \cite{FLW19,JW15,LZFW21,LFWC19,WB12,ZZFL21}. However, due to the singularity caused by the fractional operators, it is preferred to use nonuniform spatial grids to better capture the singular behavior,  and  that would destroy the special matrix structure and prevent the use of previous developed fast algorithms.
	\end{remark}
	
	\begin{remark}
		For $1<\alpha^*\leq\alpha<2$ $ (\alpha^*\approx 1.3543)$ and $\frac{\tau}{h^\alpha}<\min\{\frac{4}{\gamma},\frac{4}{1-\gamma}\}$, it can be proved that the coefficient matrix in \eqref{equ:matrix form}  is positive define if constant coefficients and uniform grids are considered.
		Thus,  the fractional CN-BCFD method \eqref{equ:sche-cn-bcfd} is uniquely solvable in such case. However,  more advanced mathematical techniques are required for the  analysis of non-uniform grids and general diffusion coefficients, which presents a significant challenge for future work.
	\end{remark}
	
	\section{A fast fractional CN-BCFD method on general nonuniform grids}\label{sec:fbcfd}
	
		\setcounter{section}{3} \setcounter{equation}{0} 

In this section, we present a fast version fractional CN-BCFD method based on Krylov subspace iterative methods for the SFDEs \eqref{model-2d:e1}. For ease of exposition, we prefer to utilize the biconjugate gradient stabilized (BiCGSTAB) method, which has faster and smoother convergence than other Krylov subspace methods. 

We remark that in the standard BiCGSTAB method, the evaluations of  matrix-vector multiplications  require $\mathcal{O}(M^2)$ computational complexity per iteration and $\mathcal{O}(M^2)$ memory for a dense and full matrix, while all other computations  require only $\mathcal{O}(M)$ computational complexity and $\mathcal{O}(M)$ memory. Therefore, it is essential to construct a  fast fractional CN-BCFD method based on the BiCGSTAB solver by developing  a fast matrix-vector multiplication mechanism and an efficient matrix storage approach.

	In this section, the SOE technique \cite{JZZZ17}, which was originally proposed to fast evaluation of the time-fractional derivative, combined with the developed fractional BCFD method will be applied to approximate the space-fractional diffusion equations. To develop a mechanism for  efficient storage of the coefficient matrix  and fast matrix-vector multiplications, we also split the left/right-sided Riemann-Liouville fractional integrals into two parts -- a local part and a history part, i.e.,
	\begin{equation}\label{equ:left-frac-inte}
		\begin{aligned}
			g^L\left(x_i,t_n\right)&=\int_a^{x_{i-1}}\frac{\left(x_i-\xi\right)^{1-\alpha}}{\Gamma(2-\alpha)}u(\xi,t_n)d\xi+\int_{x_{i-1}}^{x_i}\frac{\left(x_i-\xi\right)^{1-\alpha}}{\Gamma(2-\alpha)}u(\xi,t_n)d\xi\\
			&=:\mathcal{J}_{his,i}^{L}[u^n]+\mathcal{J}_{loc,i}^{L}[u^n],
		\end{aligned}
	\end{equation}
	\begin{equation}\label{equ:right-frac-inte}
		\begin{aligned}
			g^R\left(x_i,t_n\right)&=\int_{x_i}^{x_{i+1}}\frac{\left(\xi-x_i\right)^{1-\alpha}}{\Gamma(2-\alpha)}u(\xi,t_n)d\xi+\int_{x_{i+1}}^b\frac{\left(\xi-x_i\right)^{1-\alpha}}{\Gamma(2-\alpha)}u(\xi,t_n)d\xi\\
			&=:\mathcal{J}_{loc,i}^{R}[u^n]+\mathcal{J}_{his,i}^{R}[u^n],
		\end{aligned}
	\end{equation}
	for $1\leq i\leq M$.
	
	Without loss of generality, below we only pay attention to the derivation of the approximation formula for the left-sided Riemann-Liouville fractional integral by the SOE technique, while it is analogous for the right-sided Riemann-Liouville fractional integral.    First, for the case $i=1$, the left-sided Riemann-Liouville fractional integral can be calculated directly as in Section \ref{sec:dbcfd} since it only  has local part
	\begin{equation}\label{equ:fast-inte-1}
		g^L\left(x_1,t_n\right)\approx q^L_{1,1}u_1^n+q_{1,2}^Lu_2^n=:\mu_{1}^L u_{1}^n+\nu_1^L  u_2^n.
	\end{equation}

For $2\leq i\leq M$, as the local part $\mathcal{J}_{loc,i}^{L}[u^n]$ contributes few memory and computational cost compared with the history part $\mathcal{J}_{his,i}^{L}[u^n]$, we calculate the local part directly using \eqref{equ:interpolation} as
	\begin{equation}\label{equ:local}
		\mathcal{J}_{loc,i}^{L}[u^n]\approx\int_{x_{i-1}}^{x_i}\frac{\left(x_i-\xi\right)^{1-\alpha}}{\Gamma(2-\alpha)} \mathcal{L}_h[u^n](\xi)d\xi
		=\mu_{i}^L u_{i-1}^n+\nu_i^L u_i^n,
	\end{equation}
	where 
	\begin{equation}\label{equ:fast-inte-2}
		\mu_i^{L}=\left[\frac{1}{\Gamma(3-\alpha)}-\frac{1}{\Gamma(4-\alpha)}\right]h_{i-1/2}^{2-\alpha},\quad\nu_i^L=\frac{h_{i-1/2}^{2-\alpha}}{\Gamma(4-\alpha)},\quad 2\leq i\leq M.
	\end{equation}
Then, all the remains are to approximate the integral on the interval $\left[a, x_{i-1}\right]$ in \eqref{equ:left-frac-inte} efficiently and accurately. 	To this aim, we revisit the SOE approximation for the weak singularity kernel function $x^{1-\alpha}$ below. 
	\begin{lemma}[\cite{JZZZ17}]\label{equ:lemmasoe}
		For  given $\alpha\in(1,2)$, an absolute tolerance error $\epsilon$,  a cut-off restriction $\Delta x>0$ and a given position $X>0$, there exists a positive integer $N_{exp}$, positive quadrature points $\left\lbrace \lambda_s\right\rbrace_{s=1}^{N_{exp}} $ and corresponding positive weights $\left\lbrace \theta_s\right\rbrace_{s=1}^{N_{exp}}$ satisfying
		\begin{equation}\label{equ:kernal-soe-x}
			\Big| x^{1-\alpha} -\sum_{s=1}^{N_{exp}}\theta_se^{-\lambda_sx}\Big|\leq \epsilon,\quad\forall x\in[\Delta x,X],
		\end{equation}
		where the number of exponentials satisfies
\[
		N_{exp}=\mathcal{O}\Big(\log\frac{1}{\epsilon}\big(\log\log\frac{1}{\epsilon}+\log\frac{X}{\Delta x}\big)
		        +\log\frac{1}{\epsilon} \big(\log\log\frac{1}{\epsilon}+\log\frac{1}{\Delta x}\big)\Big).
		\]
	\end{lemma}
	
	Motivated by the above lemma, 
	we replace the convolution kernel $x^{1-\alpha}$ by its SOE approximation in \eqref{equ:kernal-soe-x} and $u(\xi,t_n)$ by its linear interpolation  $\mathcal{L}_h[u^n](\xi)$ in \eqref{equ:interpolation}, then  the history part in \eqref{equ:left-frac-inte} is approximated as follows:
	\begin{equation}\label{equ:history}
		\begin{aligned}
		\mathcal{J}_{his,i}^{L}[u^n]&\approx\frac{1}{\Gamma(2-\alpha)}\int_a^{x_{i-1}}\mathcal{L}_h[u^n](\xi) \sum_{s=1}^{N_{exp}}\theta_se^{-\lambda_s\left(x_i-\xi\right)}d\xi\\ &=\frac{1}{\Gamma(2-\alpha)}\sum_{s=1}^{N_{exp}}\theta_s\mathcal{S}_{i,s}^{L,n},
		\end{aligned}
	\end{equation}
	for  $2\leq i\leq M$,  where 
	\begin{equation}\label{equ:his:e2}
		\mathcal{S}_{i,s}^{L,n}:=\int_a^{x_{i-1}}\mathcal{L}_h[u^n](\xi)e^{-\lambda_s\left(x_i-\xi\right)}d\xi,\quad 2\leq i\leq M.
	\end{equation}
	A simple calculation shows that  the integral  in \eqref{equ:his:e2} can be computed recursively as
	\begin{equation}\label{equ:recur-rela-l2}
		\mathcal{S}_{2,s}^{L,n}=\int_a^{x_1}e^{-\lambda_s\left(x_2-\xi\right)}\mathcal{L}_h[u^n]\left(\xi\right)d\xi=\rho^L_{2,s}u_1^n+\sigma^L_{2,s}u_2^n, \\
	\end{equation}	
	\begin{equation}\label{equ:recur-rela}
		\begin{aligned}
			\mathcal{S}_{i,s}^{L,n}&=\int_a^{x_{i-2}}\mathcal{L}_h[u^n](\xi)e^{-\lambda_s\left(x_i-\xi\right)}d\xi+\int_{x_{i-2}}^{x_{i-1}}\mathcal{L}_h[u^n](\xi)e^{-\lambda_s\left(x_i-\xi\right)}d\xi\\
			&=e^{-\lambda_sh_{i-1/2}}\mathcal{S}_{i-1,s}^{L,n}+\int_{x_{i-2}}^{x_{i-1}}\mathcal{L}_h[u^n](\xi)e^{-\lambda_s\left(x_i-\xi\right)}d\xi\\
			&=e^{-\lambda_sh_{i-1/2}} \mathcal{S}_{i-1,s}^{L,n}+\rho^L_{i,s}u_{i-2}^n+\sigma^L_{i,s}u_{i-1}^n,\quad 3\leq i\leq M,
		\end{aligned}
	\end{equation}
	where
	\begin{equation}\left\{
		\begin{aligned}
			&\rho_{2,s}^L=\frac{2h_1+h_2}{h_1+h_2}\eta_{2,s}^L+\zeta_{2,s}^L,\quad \sigma_{2,s}^L=-\frac{h_1}{h_1+h_2}\eta_{2,s}^L,\\
			&\rho_{i,s}^L=\eta_{i,s}^L, \quad \sigma_{i,s}^L=\zeta_{i,s}^L,\quad 3\leq i\leq M,
		\end{aligned}\right.
	\end{equation}
	with
	\begin{equation*}\left\{
		\begin{aligned}
			\eta_{i,s}^L&=\frac{e^{-\lambda_sh_{i-1/2}}}{\lambda_s}\left(-e^{-\lambda_sh_{i-3/2}}+\frac{1-e^{-\lambda_sh_{i-3/2}}}{\lambda_sh_{i-3/2}}\right),\quad 2\leq i\leq M,\\
			\zeta_{i,s}^L&=\frac{e^{-\lambda_sh_{i-1/2}}}{\lambda_s}\left(1-\frac{1-e^{-\lambda_sh_{i-3/2}}}{\lambda_sh_{i-3/2}}\right),\quad 2\leq i\leq M.
		\end{aligned}\right.
	\end{equation*}
	
	\begin{remark}	Note that at each time step, we only need $\mathcal{O}(1)$ work to compute $\mathcal{S}_{i,s}^{L,n}$, as $\mathcal{S}_{i-1,s}^{L,n}$ is already known at that point. Therefore,  the evaluation of the convolution \eqref{equ:his:e2} can be accelerated via the recurrence relation \eqref{equ:recur-rela-l2}--\eqref{equ:recur-rela}.
	\end{remark}

	Let $g_i^{f,L,n}$  be the fast numerical approximation to left-sided Riemann-Liouville fractional integral for $1\leq i\leq M$. Therefore, combining \eqref{equ:left-frac-inte}, \eqref{equ:fast-inte-1}, \eqref{equ:local} and \eqref{equ:history} together, we obtain a fast evaluation formula for the left-sided Riemann-Liouville fractional integral as
	\begin{equation}\label{equ:fast-inte-appr-l}
		\left\{
		\begin{aligned}
			{g}_{1}^{f,L,n}&=\mu_{1}^L u_{1}^n+\nu_1^L  u_2^n,\\
			{g}_{i}^{f,L,n}&=\frac{1}{\Gamma(2-\alpha)}\sum_{s=1}^{N_{exp}}\theta_s\mathcal{S}_{i,s}^{L,n}+\mu^L_{i}u_{i-1}^n+\nu^L_{i}u_i^n,\quad 2\leq i\leq M, 
		\end{aligned}\right.
	\end{equation}
	with $\mathcal{S}_{i,s}^{L,n}$, $2\leq i \leq M$ computed recursively by \eqref{equ:recur-rela-l2}--\eqref{equ:recur-rela} for each index $s$.
	
	Analogously, 	let   $g_i^{f,R,n}$ be the fast numerical approximation to right-sided Riemann-Liouville fractional integral for $1\leq i\leq M$.  By using the same approach as above, the fast approximation for the right-sided Riemann-Liouville fractional integral is proposed as follows
	\begin{equation}\label{equ:fast-inte-appr-r}
		\left\{
		\begin{aligned}
			{g}_{M}^{f,R,n}&= 
			\nu^R_{M}u_{M-1}^n+\mu^R_{M}u_M^n,\\
			{g}_{i}^{f,R,n}&=\frac{1}{\Gamma(2-\alpha)}\sum_{s=1}^{N_{exp}}\theta_s\mathcal{S}_{i,s}^{R,n}+\nu^R_{i}u_i^n+\mu^R_{i}u_{i+1}^n,\quad 1\leq i\leq M-1,
		\end{aligned}\right.
	\end{equation}
	where 
	\begin{equation}\label{equ:fast-inte-3}
		\begin{aligned}
			&\mu^R_{M}=q^R_{M,M}, \quad \nu^R_{M}=	q^R_{M,M-1}, \\
			&\mu_i^{R}=\Big[\frac{1}{\Gamma(3-\alpha)}-\frac{1}{\Gamma(4-\alpha)}\Big]h_{i+1/2}^{2-\alpha}, ~\nu_i^R=\frac{h_{i+1/2}^{2-\alpha}}{\Gamma(4-\alpha)}, ~ 1\leq i\leq M-1,
		\end{aligned}
	\end{equation}
	and  the integral
	$\big\lbrace \mathcal{S}_{i,s}^{R,n}\big\rbrace $ satisfies the following recursive formula  
	\begin{equation}\label{equ:recur-rela-r}	\left\{
	\begin{aligned}
		&\mathcal{S}_{M-1,s}^{R,n}=\rho^R_{M-1,s}u_{M-1}^n+\sigma^R_{M-1,s}u_M^n,\\
		&\mathcal{S}_{i,s}^{R,n}=e^{-\lambda_sh_{i+1/2}}\mathcal{S}_{i+1,s}^{R,n}+\rho^R_{i,s}u_{i+2}^n+\sigma^R_{i,s}u_{i+1}^n, \quad  M-2 \ge i \ge 1.
	\end{aligned}	\right.
	\end{equation}
	The coefficients in \eqref{equ:recur-rela-r} are defined by
	\begin{equation*}\left\{
		\begin{aligned}
			& \rho_{M-1,s}^R=-\frac{h_M}{h_{M-1}+h_M}\eta_{M-1,s}^R,\quad \sigma_{M-1,s}^R=\frac{2h_{M-1}+h_M}{h_{M-1}+h_M}\eta_{M-1,s}^R+\zeta_{M-1,s}^R, \\
			&\rho_{i,s}^R=\eta_{i,s}^R,\quad \sigma_{i,s}^R=\zeta_{i,s}^R,\quad 1\leq i\leq M-2,
		\end{aligned}\right.
	\end{equation*}
	with
	\begin{equation*}\left\{
		\begin{aligned}
			\eta_{i,s}^R&=\frac{e^{-\lambda_sh_{i+1/2}}}{\lambda_s}\left(-e^{-\lambda_sh_{i+3/2}}-\frac{e^{-\lambda_sh_{i+3/2}}-1}{\lambda_sh_{i+3/2}}\right),\quad 1\leq i\leq M-1,\\
			\zeta_{i,s}^R&=\frac{e^{-\lambda_sh_{i+1/2}}}{\lambda_s}\left(1+\frac{e^{-\lambda_sh_{i+3/2}}-1}{\lambda_sh_{i+3/2}}\right),\quad 1\leq i\leq M-1.
		\end{aligned}\right.
	\end{equation*}
	
	Now,  by combining the above approximations together  with \eqref{equ:time-semi-discrete},  we get the following fast version fractional CN-BCFD  scheme: 
	\begin{equation}\label{equ:sche-fcn-bcfd}
		\left\{
		\begin{aligned}{}
			&\delta_tu_i^{n}- \frac{1}{2}  D_x p_i^{n}=\frac{1}{2}  D_x p_i^{n-1}+f_i^{n-1/2}, \quad 1\leq i\leq M,\\
			& p^{n}_{i+1/2}=\gamma K^{L,n}_{i+1/2} \, d_x g^{f,L,n}_{i+1/2}+(1-\gamma)K^{R,n}_{i+1/2}\, d_xg^{f,R,n}_{i+1/2}, \quad 1\leq i\leq M-1,\\
			&  p^{n}_{1/2}=\phi^{n},\quad
			p^{n}_{M+1/2}=\varphi^{n}, \\
			& u_i^{0}=u^o(x_i),\quad 1\leq i\leq M,
		\end{aligned}
		\right.
	\end{equation}
	with $g_{i}^{f,L,n}$ and $g_{i}^{f,R,n}$ defined by \eqref{equ:fast-inte-appr-l} and \eqref{equ:fast-inte-appr-r}, respectively.
	
	Again, let $\bm{u}^n$ and $	\bm{F}^{n-1/2} $ be defined as in Section \ref{sec:dbcfd}, we can rewrite the fast version fractional CN-BCFD scheme \eqref{equ:sche-fcn-bcfd} into the following matrix form: 
	\begin{equation}\label{equ:matrix form:f}
		\left(\bm{I}_M-\frac{\tau}{2}\bm{B}^n\right)\bm{u}^{n}=\left(\bm{I}_M+\frac{\tau}{2}\bm{B}^{n-1}\right)\bm{u}^{n-1}+\tau\bm{F}^{n-1/2},
	\end{equation}
	where the stiffness matrix	$\bm{B}^{n}$  of order $M$ has a special matrix representation
	\begin{equation}\label{equ:matrix-a}
		\bm{B}^{n} =\gamma\bm{B}^{L,n}+(1-\gamma)\bm{B}^{R,n}, 
	\end{equation}
	with $\bm{B}^{n}\bm{u}^{n}$ corresponds to  fast discretization of $\partial_x \, p(x,t_n)$, and the entries of the matrix-vector multiplication $\bm{B}^{n}\bm{u}^{n}$ are given by
	\begin{equation}\label{equ:AL}
	\left\{
		\begin{aligned}
			&(\bm{B}^{n}\bm{u}^{n})_1 = \frac{p_{3/2}^{n} }{h_1}\\
			         & \quad=\frac{\gamma}{h_1}\big[ \mathcal{K}^{L,n}_{3/2} ({g}_{2}^{f,L,n}-{g}_{1}^{f,L,n} )\big]
			                 +\frac{(1-\gamma)}{h_1} \big[ \mathcal{K}^{R,n}_{3/2} ({g}_{2}^{f,R,n}-{g}_{1}^{f,R,n} )\big],\\
			&(\bm{B}^{n}\bm{u}^{n})_i=\frac{p^{n}_{i+1/2}-p_{i-1/2}^{n}}{h_i}\\
			      &\quad =\frac{\gamma}{h_i}\big[ \mathcal{K}^{L,n}_{i+1/2}  ({g}_{i+1}^{f,L,n}-{g}_{i}^{f,L,n} )
			                    - \mathcal{K}^{L,n}_{i-1/2}  ({g}_{i}^{f,L,n}-{g}_{i-1}^{f,L,n} )\big]\\
		                    	&\qquad
		                   +\frac{(1-\gamma)}{h_i}\big[ \mathcal{K}^{L,n}_{i+1/2}  ({g}_{i+1}^{f,R,n}-{g}_{i}^{f,R,n} )
		                   - \mathcal{K}^{R,n}_{i-1/2} ({g}_{i}^{f,R,n}-{g}_{i-1}^{f,R,n})\big],\\
			&\qquad\qquad\qquad\qquad\qquad\qquad\qquad\qquad\qquad\qquad 2\leq i\leq M-1,\\
		&	(\bm{B}^{n}\bm{u}^{n})_M=-\frac{ p_{M-1/2}^{n}}{h_M} =-\frac{\gamma}{h_M}\big[ \mathcal{K}^{L,n}_{M-1/2}  ({g}_{M}^{f,L,n}-{g}_{M-1}^{f,L,n} )\big]\\
		&\qquad\qquad\qquad\qquad \qquad\quad
			    -\frac{(1-\gamma)}{h_M}\big[ \mathcal{K}^{R,n}_{M-1/2}  ({g}_{M}^{f,R,n}-{g}_{M-1}^{f,R,n} )\big].
		\end{aligned}\right.
	\end{equation}

	Based upon the above discussions, we develop a fast version BiCGSTAB iterative method for \eqref{equ:matrix form:f}, where matrix-vector multiplications  $\bm{B}^n \bm v$ for any $\bm{v}\in\mathbb{R}^M$  can be computed in an efficient way. Basically, we have the following conclusions. 
	\begin{lemma}\label{lem:BL}
		The matrix-vector multiplication $\bm{B}^{L,n}\bm{v}$ for any $\bm{v}\in\mathbb{R}^M$ can be carried out in $\mathcal{O}\left(MN_{exp}\right)$ operations, where $N_{exp}\ll M$.
	\end{lemma}
	\proof
		Let $\bm{v}=\left[v_1, v_2,\dots, v_M\right]^\top$ denote the $M$-dimensional column vector. Then, according to equation \eqref{equ:AL}, the matrix-vector multiplication $\bm{B}^{L,n}\bm{v}$ for any $\bm{v}\in\mathbb{R}^M$ is of the form
		\begin{equation}\label{equ:mv-BL}
			\begin{aligned}
				\bm{B}^{L,n}\bm{v}&=		\bm{D}_+^{L,n} \circ \big[\tilde{g}_2^{f,L}-\tilde{g}_1^{f,L},\cdots,\tilde{g}_M^{f,L}-\tilde{g}_{M-1}^{f,L},0\big]^\top\\
				&\qquad+		\bm{D}_-^{L,n}\circ \big[0,\tilde{g}_2^{f,L}-\tilde{g}_1^{f,L},\cdots,\tilde{g}_M^{f,L}-\tilde{g}_{M-1}^{f,L}\big]^\top,
			\end{aligned}
		\end{equation}
		where $ \circ $ denotes the Hadamard (element by element) product of two vectors, and $\tilde{g}_{i}^{f,L}$ has the same mathematical expression as $g_{i}^{f,L,n}$ defined by \eqref{equ:fast-inte-appr-l} just with $u_i^n$ replaced by $v_i$. Note that this step only needs $\mathcal{O}(M)$ operations. While, the computations of $\lbrace \tilde{g}_{i}^{f,L}\rbrace $ in \eqref{equ:mv-BL}, see that of $\lbrace g_{i}^{f,L}\rbrace $ in \eqref{equ:fast-inte-appr-l} require  $\mathcal{O}(MN_{exp})$ operations in total. This is due to that in the fast version fractional CN-BCFD scheme \eqref{equ:sche-fcn-bcfd}, $\mathcal{S}_{i,s}^{L}$  can be computed by the recurrence formula \eqref{equ:recur-rela} in only  $\mathcal{O}(1)$ work   for each $i$, and thus each $ g_{i}^{f,L} $ costs only $\mathcal{O}(N_{exp})$ operations. In summary, the matrix-vector multiplication $\bm{B}^{L,n}\bm{v}$ can be evaluated in $\mathcal{O}\left(MN_{exp}\right)$ operations. 
	\proofend
	
	\begin{lemma}\label{lem:BR}
		The matrix-vector multiplication $\bm{B}^{R,n}\bm{v}$ for any $\bm{v}\in\mathbb{R}^M$ can be  carried out in $\mathcal{O}\left(MN_{exp}\right)$ operations, where $N_{exp}\ll M$.
	\end{lemma}
	\proof
		Let $\bm{v}=\left[v_1, v_2,\dots, v_M\right]^\top$ denote the $M$-dimensional column vector. It is easy to show that the matrix-vector multiplication $\bm{B}^{R,n}\bm{v}$ is of the form
		\begin{equation}\label{equ:mv-BR}
			\begin{aligned}
				\bm{B}^{R,n}\bm{v}& =	\bm{D}_+^{R,n} \circ \big[\tilde{g}_2^{f,R}-\tilde{g}_1^{f,R},\cdots,\tilde{g}_M^{f,R}-\tilde{g}_{M-1}^{f,R},0\big]^\top\\
				&\qquad+		\bm{D}_-^{R,n}\circ \big[0,\tilde{g}_2^{f,R}-\tilde{g}_1^{f,R},\cdots,\tilde{g}_M^{f,R}-\tilde{g}_{M-1}^{f,R}\big]^\top,
			\end{aligned}
		\end{equation}
		where	$\tilde{g}_{i}^{f,R}$ has the same mathematical expression as $g_{i}^{f,R,n}$ defined by \eqref{equ:fast-inte-appr-r} just with $u_i^n$ replaced by $v_i$. 
		Therefore, similar to the proof of Lemma \ref{lem:BL}, we can easily conclude that the matrix-vector multiplication $\bm{B}^{R,n}\bm{v}$ can be evaluated in $\mathcal{O}\left(MN_{exp}\right)$ operations.
	\proofend
	\begin{theorem}\label{thm:fast}
		The matrix-vector multiplication $\bm{B}^{n}\bm{v}$ for any $\bm{v}\in\mathbb{R}^M$ can be  carried out in $\mathcal{O}\left(MN_{exp}\right)$ operations, where $N_{exp}\ll M$. Moreover,  the stiffness matrix $\bm{B}^{n}$ can be stored in $\mathcal{O}\left(MN_{exp}\right)$ memory. 
	\end{theorem}
	\proof
		The first conclusion can be easily obtained  by using \eqref{equ:matrix-a}, Lemmas \ref{lem:BL}--\ref{lem:BR} that the matrix-vector multiplication $\bm{B}^{n}\bm{v}$ can be performed in $\mathcal{O}\left(MN_{exp}\right)$ operations. 
		
		Furthermore, like the stiffness matrix $\bm{A}^n$ of  the fractional CN-BCFD scheme \eqref{equ:matrix form}  which is expressed  in	\eqref{coeff-a}, the matrix $\bm{B}^{n}$ also has a complicated structure. However, through the analysis of Lemmas \ref{lem:BL}--\ref{lem:BR}, we do not need to generate it explicitly. Instead, we need only to store the vectors $\bm{D}_+^{L,n}$, $\bm{D}_-^{L,n}$, $\bm{D}_+^{R,n} $ and  $\bm{D}_-^{R,n}$ defined by \eqref{coeff:e1}--\eqref{coeff:e2}, and $\lbrace \mu_i^L\rbrace $, $\lbrace \nu_i^L\rbrace $,  $\lbrace \mu_i^R\rbrace $ and $\lbrace \nu_i^R\rbrace $ defined by \eqref {equ:fast-inte-1}, \eqref{equ:fast-inte-2} and \eqref{equ:fast-inte-3}. These all require $\mathcal{O}\left(M\right)$ memory. In addition,  $\theta_s$ and $\lambda_s$ for $1\leq s\leq N_{exp}$ needs only $\mathcal{O}\left(N_{exp}\right)$ memory. Finally, the coefficients $\rho_{i,s}^L$, $\sigma_{i,s}^L$, $\rho_{i,s}^R$,  $\sigma_{i,s}^R$  for $1\leq i\leq M$, $1\leq s\leq N_{exp}$ require $\mathcal{O}\left(MN_{exp}\right)$ memory. Therefore, the total memory requirement for the stiffness matrix $\bm{B}^{n}$ is of order  $\mathcal{O}\left(M\right)+\mathcal{O}\left(N_{exp}\right)+\mathcal{O}\left(MN_{exp}\right)=\mathcal{O}\left(MN_{exp}\right)$.
	\proofend
	
	\begin{remark}	It is worth noting that in the fast version fractional CN-BCFD method, we do not need to form the matrix $\bm{B}^n$ explicitly, but rather calculate the  matrix-vector multiplication $\bm{B}^n\bm{v}$, which is different from that of Section \ref{sec:dbcfd}. The following Algorithm \ref{alg:mv1} is developed to fast carry out the matrix-vector multiplication $\bm{B}^n\bm{v}$  which appears in the BiCGSTAB method  for given vector $\bm{v}$, this yields the fast version BiCGSTAB algorithm (fBiCGSTAB) for the fractional CN-BCFD method \eqref{equ:sche-fcn-bcfd}. Note that $N_{exp}$ is usually of order $\mathcal{O}(\log^2M)$ as described in Lemma \ref{equ:lemmasoe}. Consequently, the total computational work is reduced from $\mathcal{O}(M^2)$ to $\mathcal{O}(M\log^2M)$ for each iteration of the BiCGSTAB method per time level. Meanwhile, the total memory requirement is also reduced from $\mathcal{O}(M^2)$ to $\mathcal{O}(M\log^2M)$. Thus, the fast version fractional CN-BCFD method would significantly improve the computational efficiency. Moreover, the nonuniform	spatial grids shall also improve the computational accuracy of the method.
	\end{remark}	
	\begin{algorithm}
		\caption{Fast matrix-vector multiplication $\bm{B}^n\bm{v}$ for any vector $\bm{v}$}\label{alg:mv1}
		\begin{algorithmic}[1]
			\State Compute the weights and nodes $\{\theta_s,\lambda_s\}$ of the SOE approximation using \cite{JZZZ17}
			\For {$i=2: length(\bm{v})$}
			\State Compute $\mathcal{S}_{i,s}^{L}$ by the recursive formula \eqref{equ:recur-rela-l2}--\eqref{equ:recur-rela}
			\EndFor
			\State Compute the fast approximation $\{  {g}_{i}^{f,L}\}$ by formula \eqref{equ:fast-inte-appr-l}
			\State Realize the matrix-vector multiplication $\bm{B}^{L,n}\bm{v}$ via formula \eqref{equ:mv-BL} 
			\For {$i=length(\bm{v})-1:-1:1$}
			\State Compute $\mathcal{S}_{i,s}^{R}$ by the recursive formula \eqref{equ:recur-rela-r}
			\EndFor
			\State Compute the fast approximation $\{{g}_{i}^{f,R}\}$ by formula \eqref{equ:fast-inte-appr-r}
			\State Realize matrix-vector multiplication $\bm{B}^{R,n}\bm{v}$ via formula \eqref{equ:mv-BR} 
			\State Compute $\bm{B}^n\bm{v}$ via \eqref{equ:matrix-a} that $\bm{B}^n\bm{v}=\gamma\bm{B}^{L,n}\bm{v}+(1-\gamma)\bm{B}^{R,n}\bm{v}$.
		\end{algorithmic}
	\end{algorithm}
	
	\section{Numerical experiments} \label{sec:num}
	\setcounter{section}{4} \setcounter{equation}{0} 
	
	In this section, we present some numerical examples to illustrate the convergence of the newly developed  fractional CN-BCFD  method \eqref{equ:sche-cn-bcfd} and the fast  fractional CN-BCFD  method \eqref{equ:sche-fcn-bcfd} on  nonuniform spatial grids. In addition, we also investigate the performance of the fast version fractional BCFD scheme using the developed fast BiCGSTAB iterative algorithm, in which  fast matrix-vector multiplications are carried out via Algorithm \ref{alg:mv1}. Experiments will be terminated once either the relative residual error is less than $ 10^{-10}$ or the number of iterations exceeds the maximum number of outer iterations (Here, the maximum number of outer iterations is set to be the number of spatial grids $M$). Moreover,  the tolerance error in the SOE approximation is set as $\epsilon=10^{-10}$.
	
	All the numerical experiments are performed in Matlab R2019b on a laptop with the configuration: 11th Gen Intel(R) Core (TM) i7-11700 @ 2.50GHz 2.50 GHz and 16.00 GB RAM. 
	
	\subsection{Numerical results for one-dimensional SFDEs}
	\begin{example}\label{exam:1d-1}
		For the first example, we aim to test the convergence orders of the fractional CN-BCFD scheme \eqref{equ:sche-cn-bcfd} and the fast fractional CN-BCFD scheme \eqref{equ:sche-fcn-bcfd} for enough smooth solutions. We take the	spatial interval as $[a,b]=[0,2]$ and the time interval as $[0,T]=[0,1]$.  Let $K^{L}=t x^\alpha$ and $K^{R}=t (2-x)^\alpha$ such that the true solution and its flux of  model \eqref{model-1d} are given by
		
		\begin{equation*}
				u(x,t)=e^t x^4(2-x)^4,\quad
				p(x,t)=\gamma te^t\varpi(x)-(1-\gamma)te^t\varpi(2-x),
		\end{equation*}
	\end{example}
	with $\varpi(x)=\frac{384x^{5}}{\Gamma(6-\alpha)}-\frac{3840x^{6}}{\Gamma(7-\alpha)}+\frac{17280x^{7}}{\Gamma(8-\alpha)}-\frac{40320x^{8}}{\Gamma(9-\alpha)}+\frac{40320x^{9}}{\Gamma(10-\alpha)}$.
	
	The nonuniform spatial grids used in this example are generated as follows. First, we construct a uniform partition of $[0,2]$ by $\tilde{x}_{i+1/2}=ih$, $0\leq i\leq M$ with equal grid size $h=2/M$. Then, by a small random perturbation of the grid size using the Matlab inline code, we define the nonuniform grid points as follows:
	\begin{equation}\label{equ:nonuniform-mesh}
		x_{i+1/2}=\begin{cases}0, & i=0,\\ 
			\tilde{x}_{i+1/2}+h\xi (\lambda_i-1/2), & i=1,\cdots,M-1, \\ 
			2, & i=M,\end{cases}
	\end{equation}
	where $\lambda_i$ represents a random number between $0$ and $1$, and $0\leq\xi\leq 1$ is a mesh parameter used to adjust the nonuniform grids. Obviously, when $\xi=0$, the grids generated by the above procedure \eqref{equ:nonuniform-mesh} are uniform. With the increasing of $\xi$, the nonuniformity of the spatial grids gradually increases and reaches its maximum at $\xi=1$. It is worth noting that, for each different $M$, the nonuniform grids are generated randomly.
	
	\begin{table}[!htbp]
		\centering 
		\caption{Errors and convergence orders for the fractional CN-BCFD scheme and its fast version with fixed $\alpha=1.8$, $\gamma=0.5$ } \label{tab:1d-ex1}
		\begin{tabular}{ c| c| c c |c c }
			\toprule
			Method      & $M$   & Error-u    & Cov.   & Error-p    & Cov. \\	
			\midrule
			& $2^5$ & 1.8523e-02 & ---    & 7.8121e-03 & ---     \\ 
			fractional      & $2^6$ & 4.5784e-03 & 2.0164 & 1.9189e-03 & 2.0255 \\ 
			CN-BCFD	        & $2^7$ & 1.3785e-03 & 1.7317 & 4.9696e-04 & 1.9491 \\ 
			& $2^8$ & 3.5203e-04 & 1.9694 & 1.1973e-04 & 2.0533 \\
			& $2^9$ & 1.0726e-04 & 1.7146 & 2.9218e-05 & 2.0349 \\
			\midrule
			& $2^5$ & 1.8523e-02 & ---    & 7.8121e-03 & ---     \\ 
			fast fractional	& $2^6$ & 4.5784e-03 & 2.0164 & 1.9189e-03 & 2.0255 \\ 
			CN-BCFD	    & $2^7$ & 1.3785e-03 & 1.7317 & 4.9696e-04 & 1.9491 \\ 
			& $2^8$ & 3.5203e-04 & 1.9694 & 1.1973e-04 & 2.0533 \\
			& $2^9$ & 1.0726e-04 & 1.7146 & 2.9214e-05 & 2.0350 \\
			\bottomrule
		\end{tabular}
	\end{table}
	
	\begin{figure}[!htbp]
		\centering
		\subfigure[fractional CN-BCFD]{\includegraphics[width=0.44\textwidth]{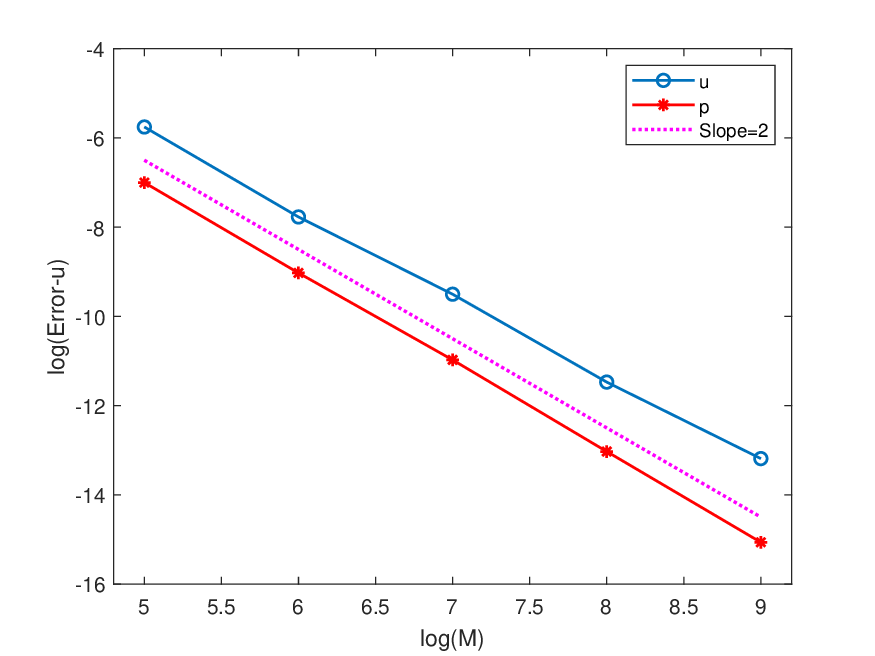}}
		\subfigure[fast fractional CN-BCFD]{\includegraphics[width=0.44\textwidth]{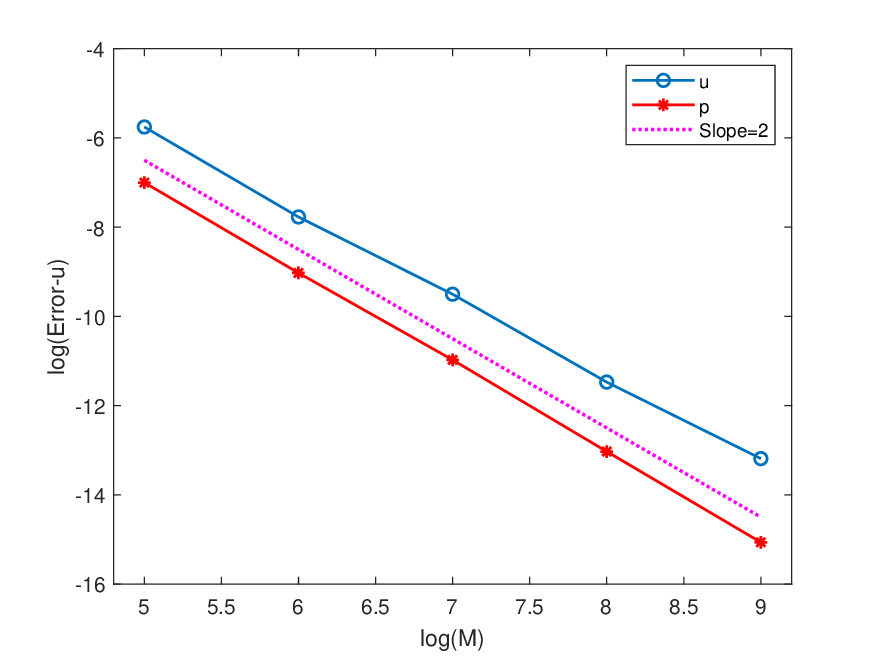}}
		\caption{Errors and convergence orders  for the fractional CN-BCFD scheme and its fast version with fixed $\alpha=1.8$, $\gamma=0.5$}
		\label{fig_ex1}
	\end{figure}	
	In this simulation, we always take $N=M$  to test the errors of the primal variable $u$ (denoted as 'Error-u') and its flux $p$ (denoted as 'Error-p') as well as their convergence orders (Cov.). In Table  \ref{tab:1d-ex1}, we list the discrete maximum-norm errors and convergence orders with the chosen parameter $\xi=1/3$ in \eqref{equ:nonuniform-mesh}, and we also depict the convergence orders of the developed fractional CN-BCFD scheme \eqref{equ:sche-cn-bcfd} and its fast version \eqref{equ:sche-fcn-bcfd} in Figure \ref{fig_ex1}. We can observe that the fast version fractional CN-BCFD scheme \eqref{equ:sche-fcn-bcfd} generates numerical solutions with the same accuracy as the fractional CN-BCFD scheme  \eqref{equ:sche-cn-bcfd}, and the primal variable $u$ and  its flux $p$ all have second-order convergence with respect to the discrete maximum-norm for $\alpha=1.8$ and $\gamma=0.5$  (i.e., the fractional operator is symmetric). Therefore, we conclude that, when the true solution has sufficient smoothness, the newly developed schemes can achieve second-order accuracy on  nonuniform grids.
	
	\begin{example}\label{exam:1d-2}
		In this example, we take the spatial interval as $[a,b]=[0,2]$ and the time interval as $[0,T]=[0,1]$. 
		Let $K^{x,L}=t \, (5+x^\alpha)$ and $K^{x,R}=t\left(5+(2-x)^\alpha\right)$ such that the true solution and its flux of model \eqref{model-1d} are given by

	\begin{equation*}
				u(x,t)=e^{-t}x^2(2-x)^2,\quad
				p(x,t)=\gamma te^{-t}\varpi(x)-(1-\gamma)te^{-t}\varpi(2-x),
		\end{equation*}
	\end{example}
with $\varpi(x)=\frac{40x^{3-\alpha}}{\Gamma(4-\alpha)}-\frac{120x^{4-\alpha}}{\Gamma(5-\alpha)}+\frac{120x^{5-\alpha}}{\Gamma(6-\alpha)}+\frac{8x^{3}}{\Gamma(4-\alpha)}-\frac{24x^{4}}{\Gamma(5-\alpha)}+\frac{24x^{5}}{\Gamma(6-\alpha)}$.
	
\begin{figure}[!ht]
	\vspace{-0.4cm}
	\centering
	\subfigure[$\kappa=1$]{\includegraphics[width=0.44\textwidth]{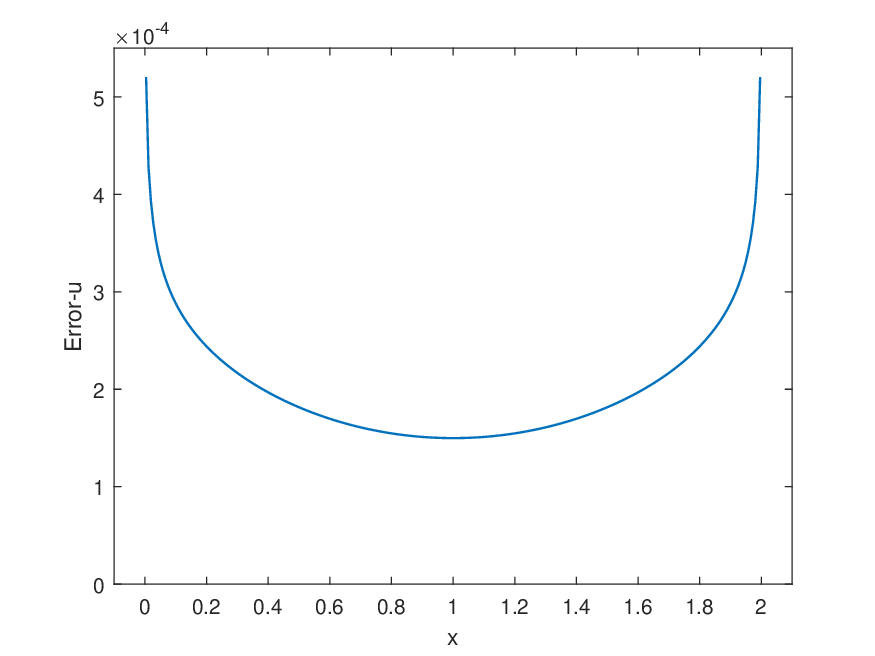}}
	\subfigure[$\kappa=1.5$]{\includegraphics[width=0.44\textwidth]{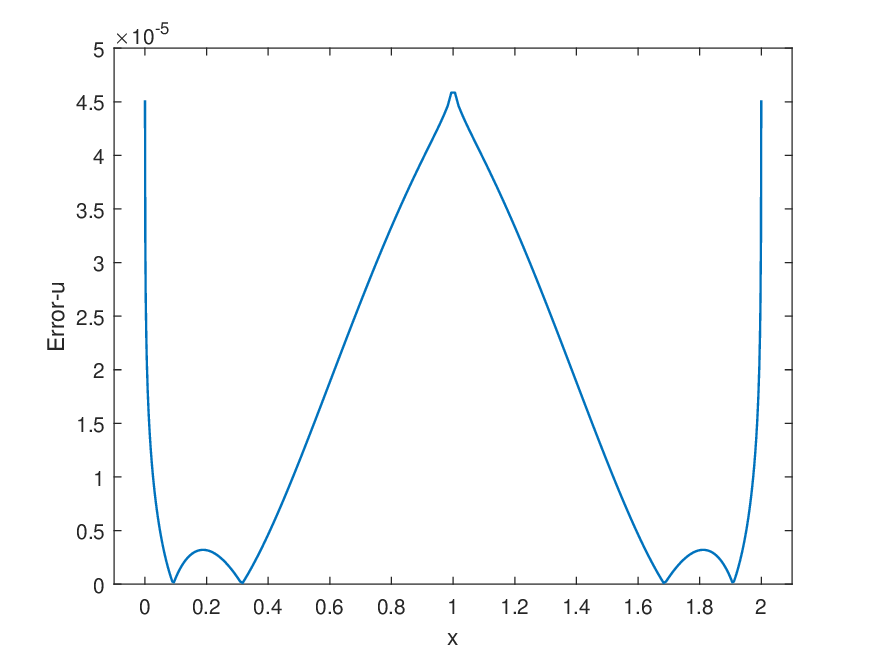}}
	\subfigure[$\kappa=1$]{\includegraphics[width=0.44\textwidth]{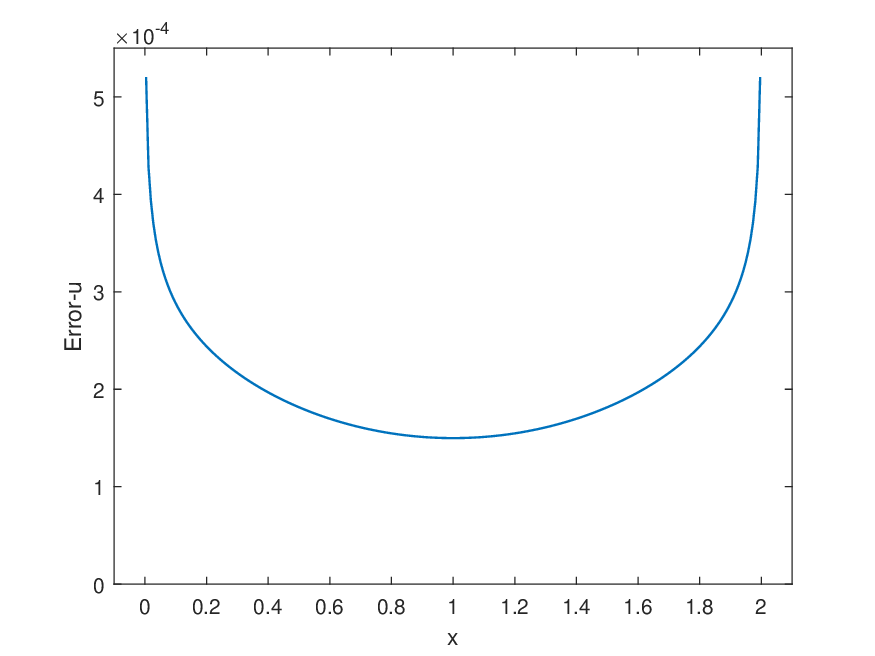}}
	\subfigure[$\kappa=1.5$]{\includegraphics[width=0.44\textwidth]{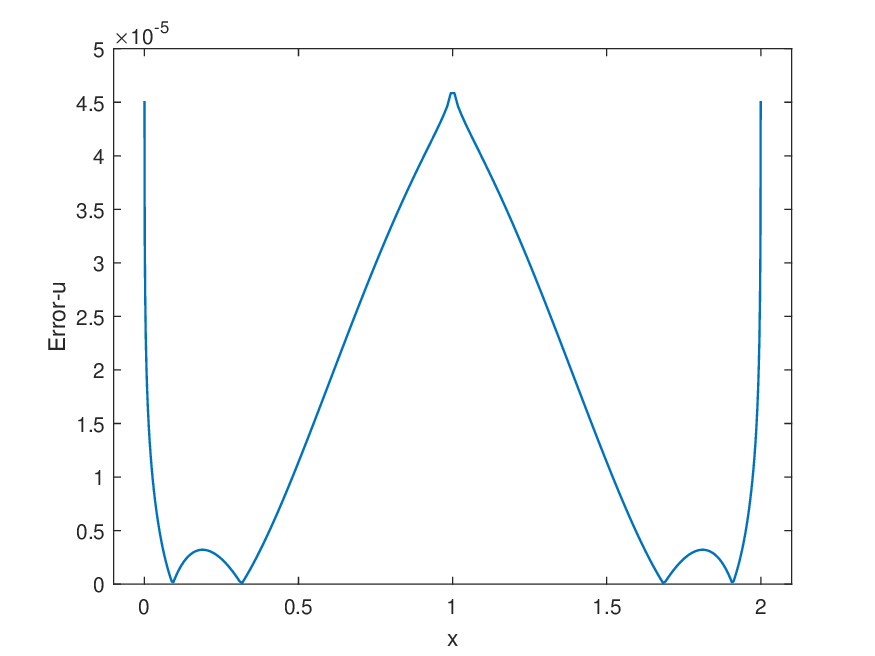}}
	\caption{Error profiles  for the fractional CN-BCFD scheme (the first row) and the fast fractional CN-BCFD scheme (the second row) with fixed $\alpha=1.5$, $\gamma=0.5$}
	\label{fig_ex2}
\end{figure}

First, we choose $M=2^8$, $N=2^{12}$ and $\alpha=1.5$, $\gamma=0.5$ to test the errors between the analytical solutions and their numerical approximations under uniform grids (i.e., $\kappa=1$ in \eqref{e2:nonuniform-mesh}). The pointwise errors of $u$ at final time $T=1$ are shown in Figure \ref{fig_ex2} (a) and (c) respectively for the fractional CN-BCFD scheme \eqref{equ:sche-cn-bcfd} and its fast version \eqref{equ:sche-fcn-bcfd}. It reveals that the error mainly comes from the neighbor of boundary points. Thus, when weak singularity occurs in model \eqref{model-1d}, the developed fractional BCFD methods with or without the fast SOE approximation may not work very well as expected and may also lose accuracy. To improve the convergence, one can employ the locally refined grids in the spatial discretization instead of uniform grids to reduce the computational errors. For example, we consider the following graded grids generated as below:
	\begin{equation}\label{e2:nonuniform-mesh}
		x_{i+1/2}=\left\{
		\begin{aligned}
			&a+\gamma(b-a) \Big(\frac{i}{ \lfloor\gamma M\rfloor } \Big)^{\kappa},    & 0\leqslant i\leqslant \lfloor\gamma M\rfloor,\\
			&b-(1-\gamma)(b-a) \Big(  \frac{M-i}{M-\lfloor\gamma M\rfloor}\Big)^{\kappa},   &\lfloor\gamma M\rfloor < i\leqslant M,
		\end{aligned}\right.
	\end{equation}
	where $\lfloor\gamma M\rfloor$ means the largest integer that less than or equal to $\gamma M$, and $\kappa$ is a user-defined graded parameter. It can be seen from Figure \ref{fig_ex2} (b) and (d) that the errors are significantly reduced when using the graded grids (i.e., $\kappa=1.5$) compared with the uniform case  (i.e., $\kappa=1$).

\begin{figure}[!htbp]
		\centering
		\subfigure[$\alpha=1.8$ ]{\includegraphics[width=0.44\textwidth]{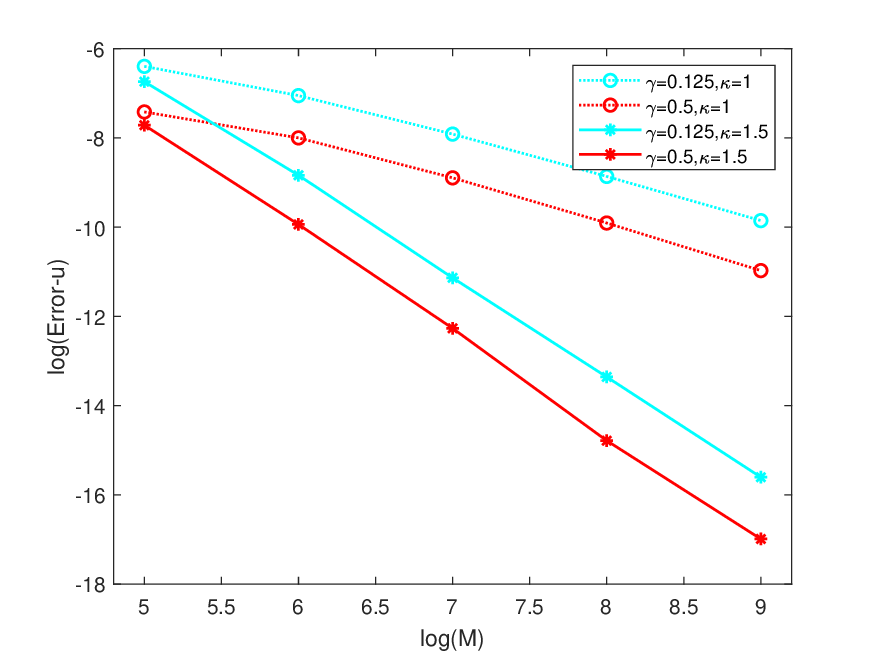}}
		\subfigure[$\alpha=1.6$ ]{\includegraphics[width=0.44\textwidth]{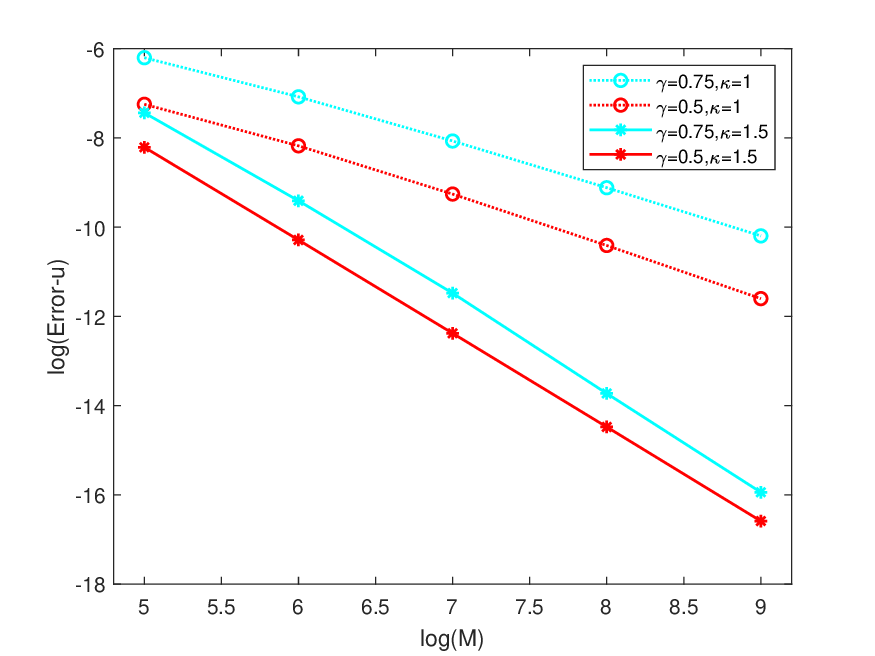}}
		\subfigure[$\alpha=1.5$ ]{\includegraphics[width=0.44\textwidth]{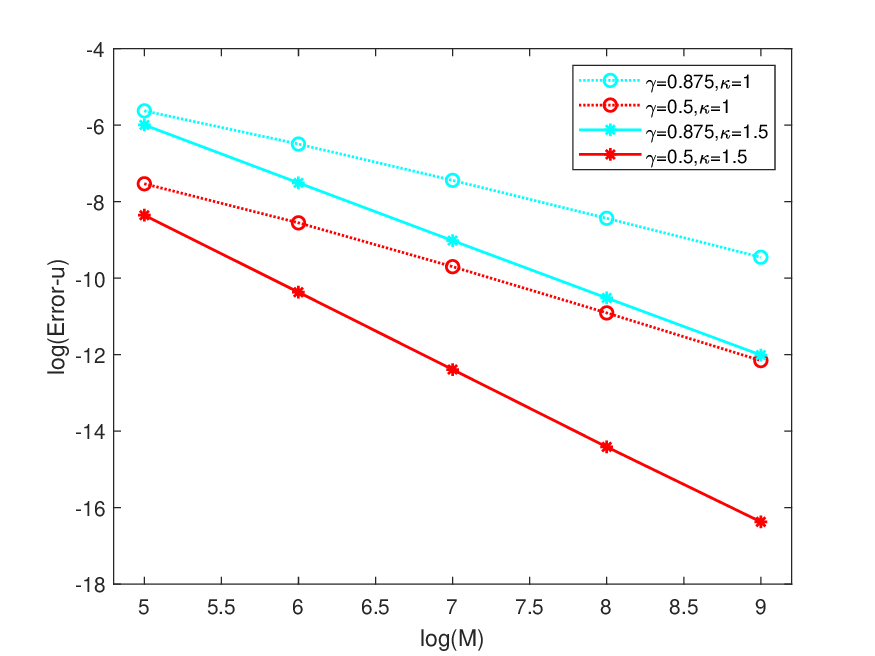}}
		\subfigure[$\alpha=1.3$ ]{\includegraphics[width=0.44\textwidth]{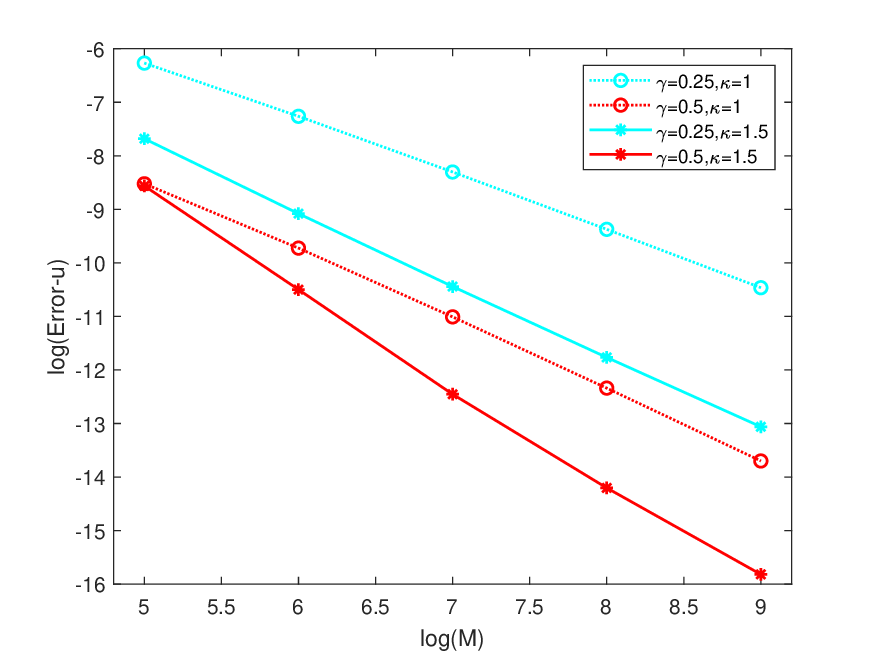}}
		\caption{{Errors of $u$ for the fast fractional CN-BCFD scheme with different fractional order $\alpha$, weighted parameter $\gamma$ and mesh parameter $\kappa$}}		\label{fig_ex21}
					\vspace{-0.3cm}
	\end{figure}

\begin{figure}[!ht]
		\centering
		\subfigure[$\alpha=1.8$]{\includegraphics[width=0.45\textwidth]{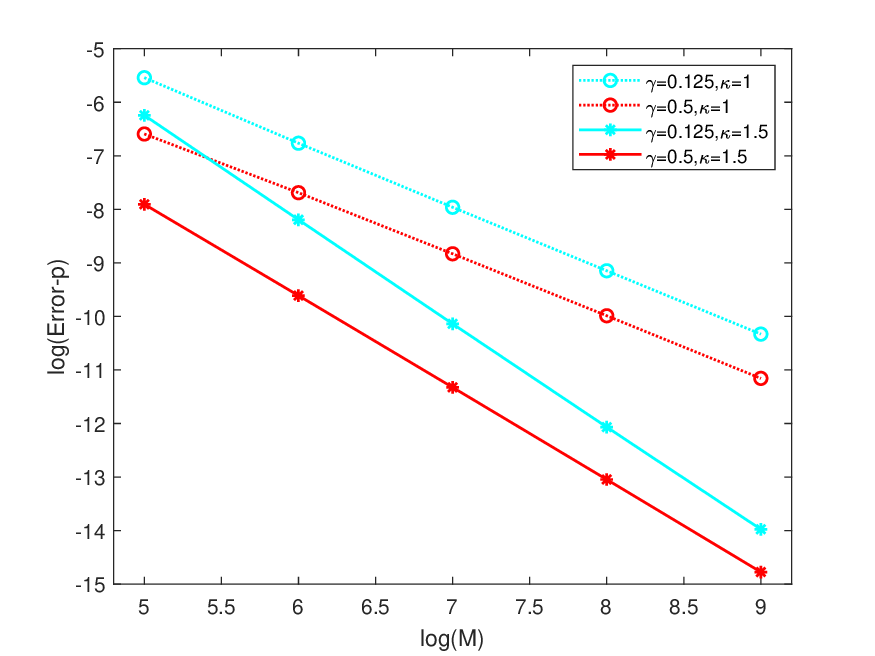}}
		\subfigure[$\alpha=1.6$]{\includegraphics[width=0.45\textwidth]{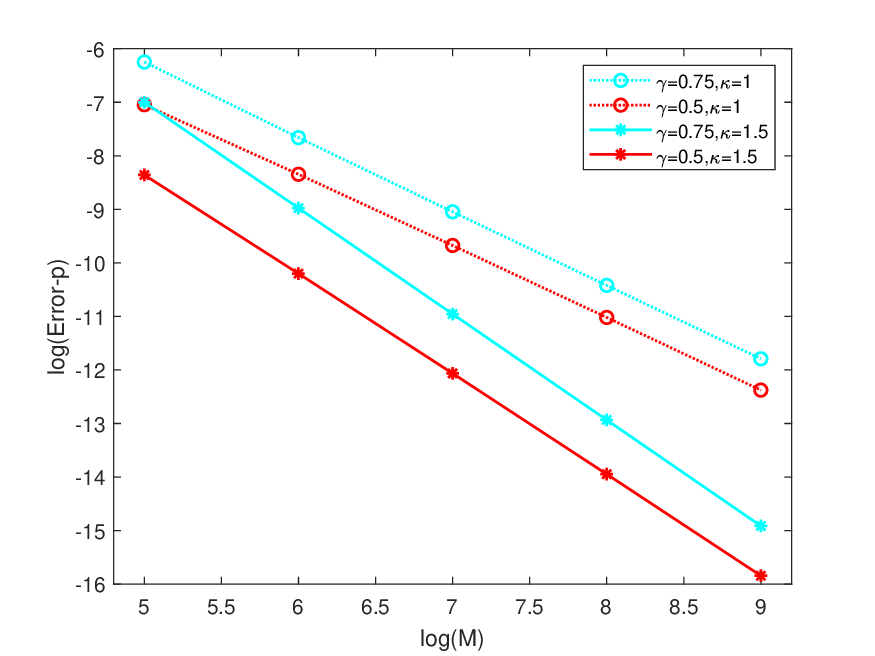}}
		\subfigure[$\alpha=1.5$]{\includegraphics[width=0.45\textwidth]{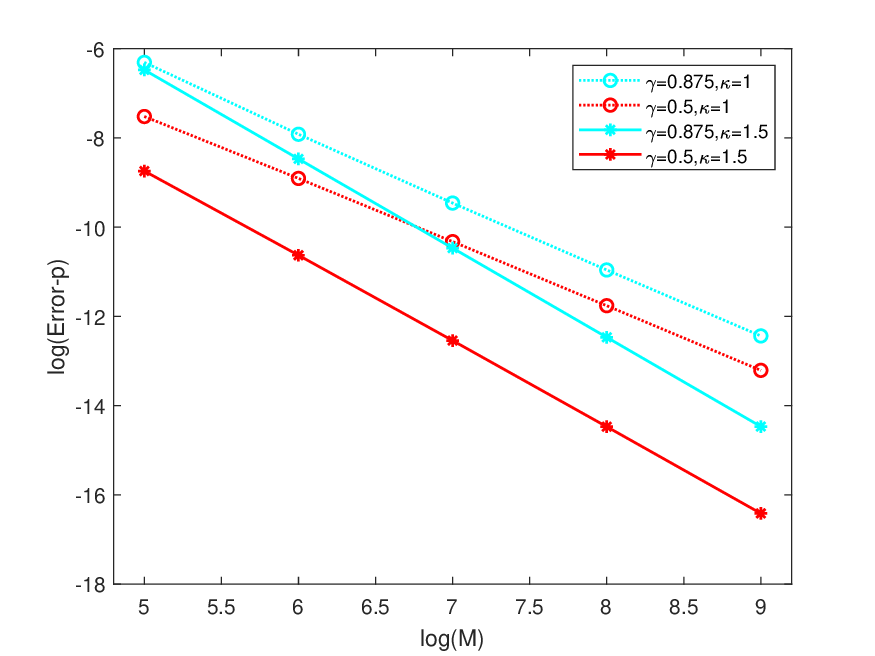}}
		\subfigure[$\alpha=1.3$]{\includegraphics[width=0.45\textwidth]{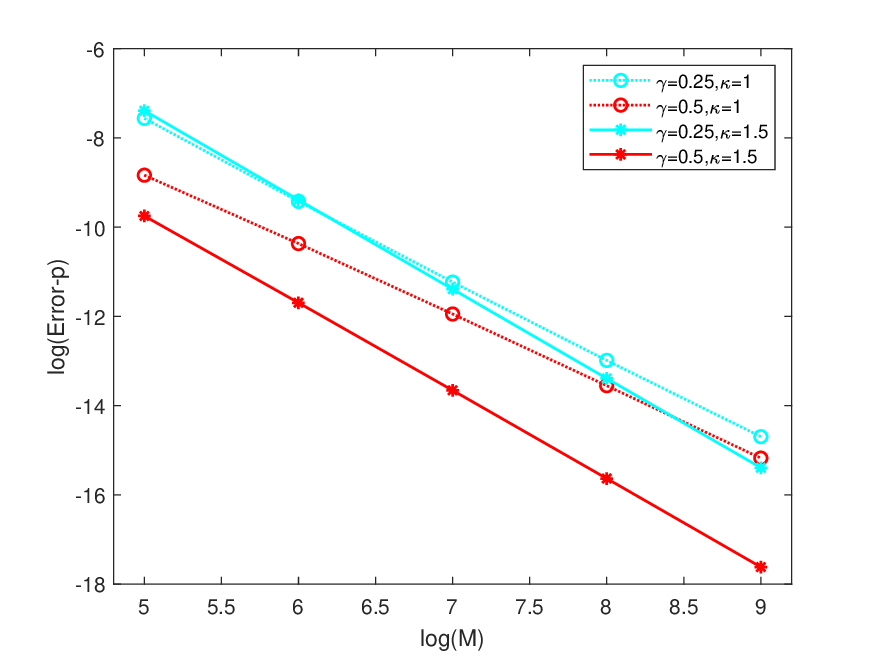}}
		\caption{Errors of $p$ for the fast fractional CN-BCFD scheme with different fractional order $\alpha$, weighted parameter $\gamma$ and mesh parameter $\kappa$}
		\label{fig_ex22}
	\end{figure}

Next, we proceed to investigate the influence of  the graded mesh parameter $\kappa$ on convergence. We choose $N = 2^{12}$ large enough to observe the spatial convergence order for Example \ref{exam:1d-2}. 	Numerical results  are listed in Tables \ref{tab:1d-ex2-u}--\ref{tab:1d-ex2-p} with respect to different choices of mesh parameters $\kappa$ for fixed $\alpha=1.5$, $\gamma=0.5$.  We can observe that the fast version fractional CN-BCFD scheme \eqref{equ:sche-fcn-bcfd} still generates the similar accurate numerical solutions as the fractional CN-BCFD scheme  \eqref{equ:sche-cn-bcfd}. Furthermore, it can be seen that for the solution that not smooth enough, the primal variable $u$ and its flux $p$ both have a poor accuracy compared with Example \ref{exam:1d-1} on the uniform grids, see Tables \ref{tab:1d-ex1} and \ref{tab:1d-ex2-u}--\ref{tab:1d-ex2-p} with $\kappa=1$. In contrast, the two schemes can achieve second-order accuracy by properly adopting suitable nonuniform spatial grids, such as the mesh parameter $\kappa=1.5$.  However, it is worth noting that the large mesh parameter may also increase the errors. For instance, Table \ref{tab:1d-ex2-u} shows that the errors generated via $\kappa=2$ is larger than that of $\kappa=1.5$. Therefore, how to select appropriate mesh parameter to achieve optimal error results is still a challenging and interesting work.  Finally,  we depict the errors of primal variable $u$ and its flux $p$ with different fractional order $\alpha$, weighted parameter $\gamma$ and mesh parameter $\kappa$ in Figures \ref{fig_ex21}--\ref{fig_ex22}. It indicates that the usage of graded grids (the mesh parameter $\kappa$) for given $\alpha$ and $\gamma$ can significantly improve the accuracy. In addition, the results obtained for the case that the fractional operator is symmetric (i.e., $\gamma=0.5$) are usually better than other cases.  In a word, when low regularity or  even weak singularity occurs in the solution, it is better to use nonuniform grids to obtain numerical results with much better accuracy, and our developed fast CN-BCFD method is well suited for general nonuniform grids. However, the reasonable choice of nonuniform grids (e.g., \eqref{e2:nonuniform-mesh}) usually depends on some a priori regularity information of the solution to \eqref{model-1d}, and as pointed out in Ref. \cite{EHR18}  that the solution to the space-fractional diffusion equations even with smooth coefficients and source terms exhibits nonphysical singularities on the spatial boundaries.

	\begin{table}[!htbp]
			\vspace{-0.4cm}
		\centering 
		\caption{Errors and convergence orders of $u$  for the fractional CN-BCFD scheme and  its fast version with fixed $\alpha=1.5$, $\gamma=0.5$}\label{tab:1d-ex2-u}
		\begin{tabular}{c |c |c c| c c| c c}
			\toprule
			\multirow{2}{*}{Method}&\multirow{2}{*}{$M$} & \multicolumn{2}{c}{$\kappa=1$} & \multicolumn{2}{|c}{$\kappa=1.5$} & \multicolumn{2}{|c}{$\kappa=2$}\\ 
			\cmidrule(r){3-4}     \cmidrule(r){5-6}     \cmidrule(r){7-8}
			&& Error-u      & Cov.    & Error-u      & Cov.    & Error-u          & Cov.    \\	
			\midrule
			&$2^5$ & 5.4058e-03 & ---     & 3.0553e-03 & ---     & 6.5626e-03     & ---     \\ 
			fractional	&$2^6$ & 2.6623e-03 & 1.0218 & 7.5680e-04 & 2.0134 & 2.0684e-03     & 1.6658 \\ 
			CN-BCFD	&$2^7$ & 1.2029e-03 & 1.1462 & 1.8648e-04 & 2.0209 & 6.5705e-04     & 1.6544 \\ 
			&$2^8$ & 5.2022e-04 & 1.2093 & 4.5861e-05 & 2.0237 & 2.0973e-04     & 1.6474 \\
			&$2^9$ & 2.1933e-04 & 1.2460 & 1.1889e-05 & 1.9477 & 6.7128e-05     & 1.6436 \\
			\midrule
			&$2^5$ & 5.4058e-03 & ---     & 3.0553e-03 & ---     & 6.5626e-03     & ---     \\ 
			fast fractional&$2^6$ & 2.6623e-03 & 1.0218 & 7.5679e-04 & 2.0134 & 2.0684e-03     & 1.6658 \\ 
			CN-BCFD&$2^7$ & 1.2029e-03 & 1.1462 & 1.8647e-04 & 2.0209 & 6.5700e-04     & 1.6545 \\ 
			&$2^8$ & 5.2022e-04 & 1.2093 & 4.5852e-05 & 2.0239 & 2.0963e-04     & 1.6481 \\
			&$2^9$ & 2.1933e-04 & 1.2460 & 1.1792e-05 & 1.9591 & 6.6807e-05     & 1.6497 \\
			\bottomrule
		\end{tabular}
	\end{table}		
	\begin{table}[!ht]
		\centering 
		\caption{Errors and convergence orders of $p$ for  the fractional CN-BCFD scheme and  its fast version with fixed $\alpha=1.5$, $\gamma=0.5$}\label{tab:1d-ex2-p}
		\begin{tabular}{c|c |c c| c c| c c}
			\toprule
			\multirow{2}{*}{Method}&\multirow{2}{*}{$M$} & \multicolumn{2}{c}{$\kappa=1$} & \multicolumn{2}{|c}{$\kappa=1.5$} & \multicolumn{2}{|c}{$\kappa=2$}\\ 
			\cmidrule(r){3-4}     \cmidrule(r){5-6}     \cmidrule(r){7-8}
			&& Error-p      & Cov.    & Error-p      & Cov.    & Error-p          & Cov.    \\	
			\midrule
			&$2^5$ & 5.4501e-03 & --     & 2.3359e-03 & --     & 2.7983e-03     & --     \\ 
			fractional&$2^6$ & 2.0828e-03 & 1.3878 & 6.3144e-04 & 1.8873 & 7.1462e-04     & 1.9693 \\ 
			CN-BCFD&$2^7$ & 7.8032e-04 & 1.4164 & 1.6775e-04 & 1.9123 & 1.8041e-04     & 1.9859 \\ 
			&$2^8$ & 2.8837e-04 & 1.4362 & 4.4015e-05 & 1.9303 & 4.5367e-05     & 1.9916 \\
			&$2^9$ & 1.0562e-04 & 1.4490 & 1.1442e-05 & 1.9436 & 1.1381e-05     & 1.9951 \\
			\midrule
			&$2^5$ & 5.4501e-03 & --     & 2.3359e-03 & --     & 2.7983e-03     & --     \\ 
			fast fractional &$2^6$ & 2.0828e-03 & 1.3878 & 6.3144e-04 & 1.8872 & 7.1462e-04     & 1.9693 \\ 
			CN-BCFD&$2^7$ & 7.8032e-04 & 1.4164 & 1.6776e-04 & 1.9123 & 1.8042e-04     & 1.9858 \\ 
			&$2^8$ & 2.8836e-04 & 1.4362 & 4.4024e-05 & 1.9300 & 4.5371e-05     & 1.9915 \\
			&$2^9$ & 1.0561e-04 & 1.4492 & 1.1456e-05 & 1.9422 & 1.1374e-05     & 1.9961 \\
			\bottomrule
		\end{tabular}
	\end{table}
	
	\begin{example}\label{exam:1d-3}
		For the third example, we proceed to investigate the performance of the proposed fractional CN-BCFD scheme \eqref{equ:sche-cn-bcfd} and the fast fractional CN-BCFD scheme \eqref{equ:sche-fcn-bcfd} with a low regularity solution
		\[
		u(x,t)=4e^tx^{\frac{\alpha}{2}}(1-x)^{\frac{\alpha}{2}},\quad x\in[0,1],\ t \in[0,1].
		\]
		Here we choose $K^{x,L}=K^{x,R}=t\left(1+x(1-x)\right)$, and the  flux $p$ is 
		\begin{equation*}
			p(x,t)=4te^t\Gamma(\alpha+1)\cos(\alpha\pi/2)\left(1+x(1-x)\right)( x-1/2).
		\end{equation*}
	\end{example}

\begin{figure}[!htbp]
	\centering
	\subfigure[$\kappa=1$]{\includegraphics[width=0.44\textwidth]{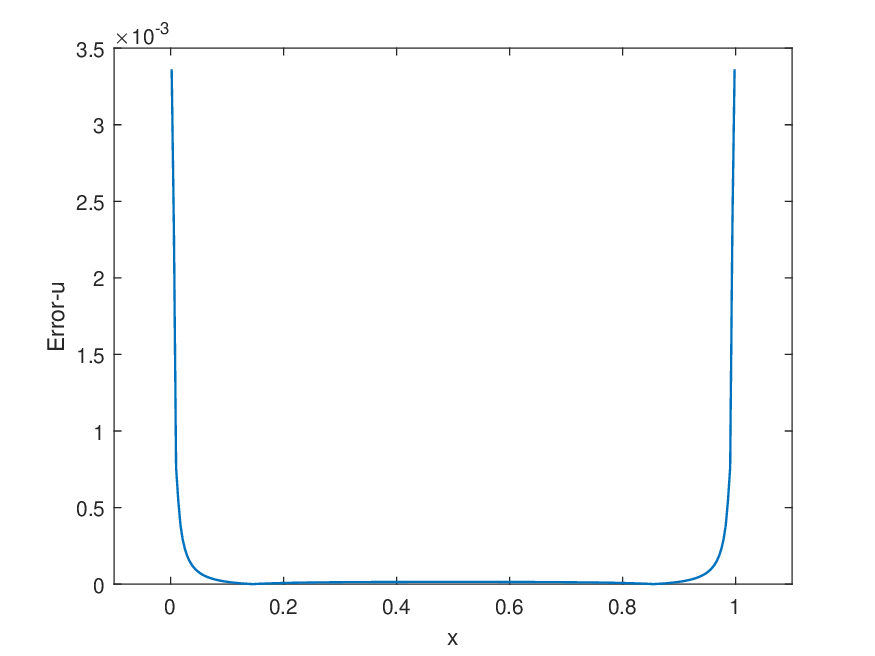}}
	\subfigure[$\kappa=1.5$]{\includegraphics[width=0.44\textwidth]{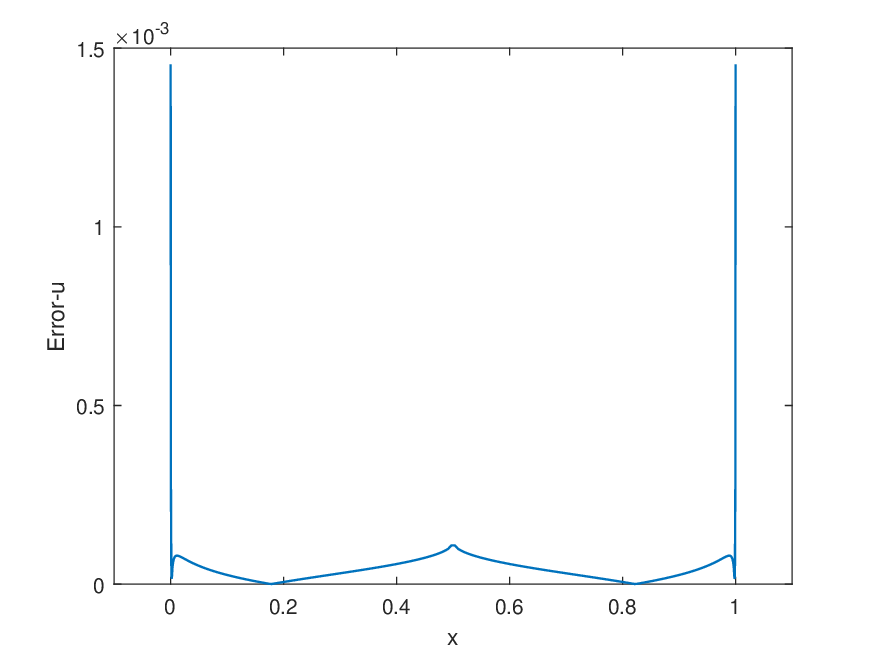}}
	\subfigure[$\kappa=1$]{\includegraphics[width=0.44\textwidth]{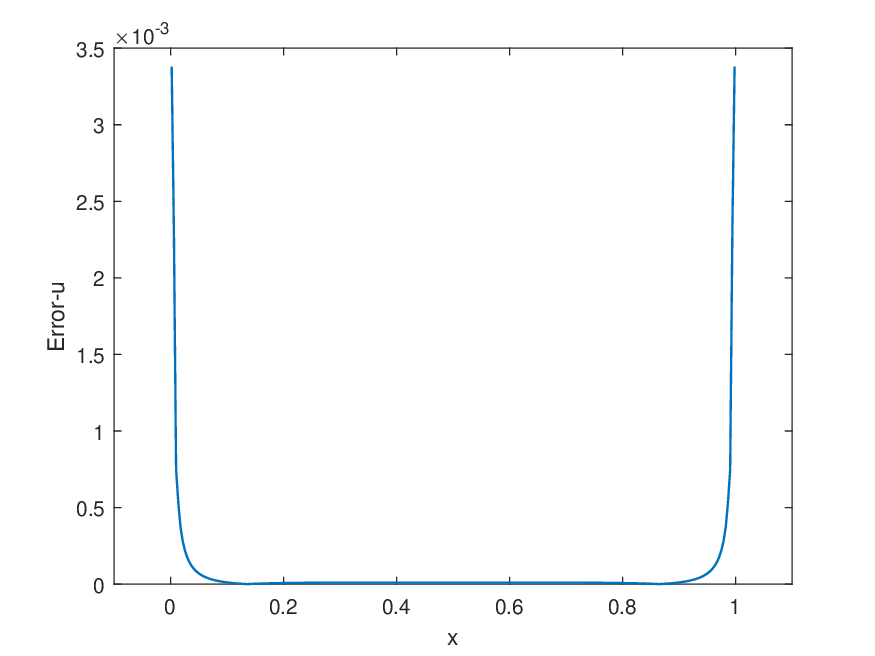}}
	\subfigure[$\kappa=1.5$]{\includegraphics[width=0.44\textwidth]{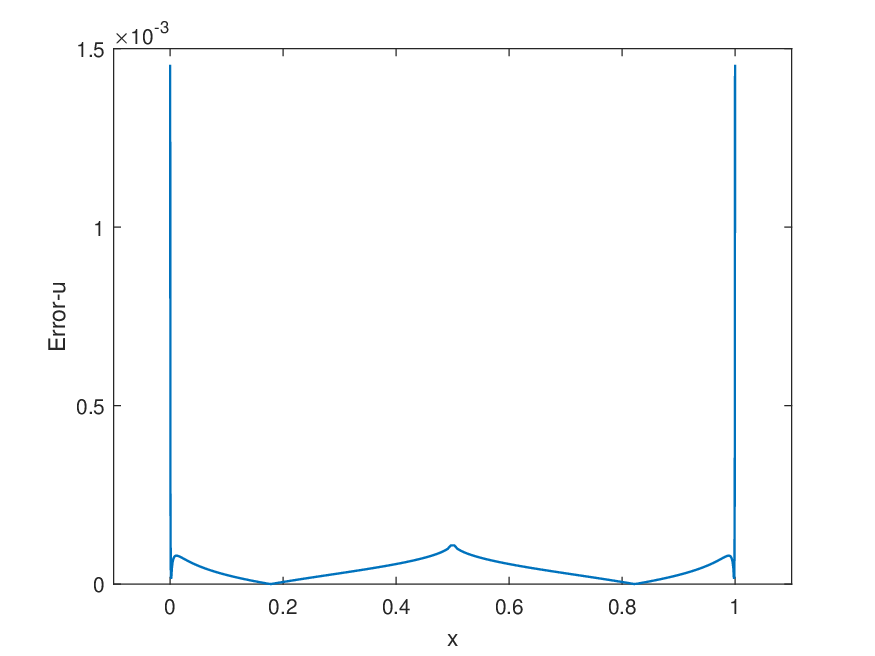}}
	\caption{Error profiles for Example \ref{exam:1d-3} at final time $T$ solved by the fractional CN-BCFD scheme (the first row) and the fast fractional CN-BCFD scheme (the second row) with fixed $\alpha=1.4$, $\gamma=0.5$}
	\label{fig_ex32}
\end{figure}		
To start with, 	 we choose $M=2^8$ and $N=2^{12}$ to test the pointwise errors between exact solutions and numerical solutions at final time $T=1$ with fixed $\alpha=1.4$ and $\gamma=0.5$. Similar as Example \ref{exam:1d-2}, we conclude from Figure \ref{fig_ex32}  that the errors are mainly from the neighbor of boundary points,  and the usage of graded grids (i.e. $\kappa=1.5$) clearly improves the accuracy. Then, we present the discrete maximum-norm errors and corresponding convergence orders in Tables \ref{tab:1d-ex3-u}--\ref{tab:1d-ex3-p}, respectively, for the primal variable $u$ and its flux $p$. It can be seen that the usage of graded grids results in significantly smaller errors than that of uniform grid (i.e. $\kappa=1$) near the boundary. Nevertheless, the convergence orders  cannot reach two for the variable $u$ despite the usage of graded grids, while the flux $p$ can achieve second-order accuracy on the graded grids. Meanwhile, the larger mesh parameter $\kappa$ will yield a reduction of accuracy. For example, the errors of $\kappa=2$ is obviously bigger than that of $\kappa=1.5$.

	\begin{table}[!htbp]
		\centering 
		\caption{Errors and convergence orders of $u$ for  the fractional CN-BCFD scheme and and its fast version with fixed $\alpha=1.4$, $\gamma=0.5$}\label{tab:1d-ex3-u}
		\begin{tabular}{c|c |c c| c c| c c}
			\toprule
			\multirow{2}{*}{Method}&\multirow{2}{*}{$M$} & \multicolumn{2}{c}{$\kappa=1$} & \multicolumn{2}{|c}{$\kappa=1.5$} & \multicolumn{2}{|c}{$\kappa=2$}\\ 
			\cmidrule(r){3-4}     \cmidrule(r){5-6}     \cmidrule(r){7-8}
			&	& Error-u      & Cov.    & Error-u      & Cov.    & Error-u          & Cov.    \\	
			\midrule
			&$2^6$ & 9.6150e-03 & --     & 1.0497e-02 & --     & 1.7792e-02 & --     \\ 
			fractional	&$2^7$ & 5.6731e-03 & 0.7611 & 3.9111e-03 & 1.4243 & 6.7802e-03 & 1.3918 \\ 
			CN-BCFD	&$2^8$ & 3.3787e-03 & 0.7477 & 1.4545e-03 & 1.4271 & 2.5782e-03 & 1.3950 \\
			&$2^9$ & 2.0309e-03 & 0.7343 & 5.4477e-04 & 1.4168 & 9.7900e-04 & 1.3970 \\
			&$2^{10}$& 1.2299e-03 & 0.7236 & 2.0697e-04 & 1.3962 & 3.7130e-04 & 1.3987 \\
			\midrule
			&	$2^6$ & 9.6152e-03 & --     & 1.0497e-02 & --     & 1.7792e-02 & --     \\ 
			fast fractional	&$2^7$ & 5.6733e-03 & 0.7611 & 3.9113e-03 & 1.4243 & 6.7805e-03 & 1.3917 \\ 
			CN-BCFD	&$2^8$ & 3.3789e-03 & 0.7476 & 1.4547e-03 & 1.4269 & 2.5796e-03 & 1.3943 \\
			&$2^9$ & 2.0315e-03 & 0.7340 & 5.4524e-04 & 1.4157 & 9.8528e-04 & 1.3885 \\
			&$2^{10}$& 1.2360e-03 & 0.7232 & 2.0791e-04 & 1.3909 & 3.5591e-04 & 1.4690 \\
			\bottomrule
		\end{tabular}
	\end{table}
	\begin{table}[!htbp]
		\centering 
		\caption{Errors and convergence orders of $p$ for  the fractional CN-BCFD scheme and its fast version with fixed $\alpha=1.4$, $\gamma=0.5$}\label{tab:1d-ex3-p}
		\begin{tabular}{c| c |c c| c c| c c}
			\toprule
			\multirow{2}{*}{Method}&\multirow{2}{*}{$M$} & \multicolumn{2}{c}{$\kappa=1$} & \multicolumn{2}{|c}{$\kappa=1.5$} & \multicolumn{2}{|c}{$\kappa=2$}\\ 
			\cmidrule(r){3-4}     \cmidrule(r){5-6}     \cmidrule(r){7-8}
			&& Error-p      & Cov.    & Error-p      & Cov.    & Error-p          & Cov.    \\	
			\midrule
			&$2^6$ & 1.4694e-04 & --     & 1.2299e-04 & --     & 3.2002e-04     & --     \\ 
			fractional&$2^7$ & 4.3279e-05 & 1.7635 & 2.8846e-05 & 2.0910 & 8.0726e-05     & 1.9871 \\ 
			CN-BCFD&$2^8$ & 1.2946e-05 & 1.7412 & 6.8784e-06 & 2.0682 & 2.0309e-05     & 1.9909 \\ 
			&$2^9$ & 3.9117e-06 & 1.7266 & 1.6597e-06 & 2.0512 & 5.0978e-06     & 1.9942 \\
			&$2^{10}$& 1.3650e-06 & 1.5189 & 4.0339e-07 & 2.0407 & 1.2767e-06     & 1.9974 \\
			\midrule
			&$2^6$ & 1.4691e-04 & --     & 1.2288e-04 & --     & 3.2002e-04     & --     \\ 
			fast fractional&$2^7$ & 4.3275e-05 & 1.7633 & 2.8830e-05 & 2.0916 & 8.0724e-05     & 1.9871 \\ 
			CN-BCFD&$2^8$ & 1.2943e-05 & 1.7414 & 6.8699e-06 & 2.0692 & 2.0322e-05     & 1.9900 \\ 
			&$2^9$ & 3.9051e-06 & 1.7287 & 1.6487e-06 & 2.0589 & 5.1767e-06     & 1.9729 \\
			&$2^{10}$& 1.4492e-06 & 1.4301 & 3.9197e-07 & 2.0725 & 9.8745e-07     & 2.3903 \\
			\bottomrule
		\end{tabular}
	\end{table}	 
	
	\begin{example}\label{exam:1d-4}
		Through this example, we mainly aim to verify the efficiency of the proposed fast fractional CN-BCFD method for the approximation of the SFDE \eqref{model-1d} with different fractional order $\alpha$ and mesh parameter $\kappa$. In this example, the diffusivity coefficients $K^{x,L}=K^{x,R}=x (2-x)$,  and the true solution is given by
		\[
		u(x,t)=4tx(2-x),\quad x\in[0,2],\ t \in[0,1].
		\]
	\end{example}
	
	In Tables \ref{tab:1d-1}--\ref{tab:1d-2},  we present the errors and  CPU times for the fractional CN-BCFD scheme \eqref{equ:sche-cn-bcfd} solved by the Gaussian elimination solver (denoted as CN-BCFD-GE), the fractional CN-BCFD scheme \eqref{equ:sche-cn-bcfd} solved by the BiCGSTAB iterative solver (denoted as CN-BCFD-BiCGSTAB) and the fast fractional CN-BCFD scheme \eqref{equ:sche-fcn-bcfd} solved by the fast version BiCGSTAB iterative solver (denoted as CN-BCFD-fBiCGSTAB) 	via Algorithm \ref{alg:mv1}. Meanwhile,  the average number of iterations (Itr.) are also shown for the iterative solver. We can reach the following observations: (i) All these methods generate almost the identical numerical solutions; (ii) The fast version CN-BCFD-fBiCGSTAB algorithm takes significantly less CPU time than the other two algorithms, especially for large-scale modeling. For example, in the case of $\alpha=1.8$, $\kappa=1$, $M=2^{11}$, the CPU time consumed by the CN-BCFD-GE algorithm is more than $46$ hours. In contrast, the fast version CN-BCFD-fBiCGSTAB algorithm only costs no more than 3 minutes; (iii) The use of nonuniform grids does not affect the computational accuracy  of the CN-BCFD-fBiCGSTAB algorithm, while compared to the other two algorithms, it is still the most efficient one. In summary, all these results show the superiority of the CN-BCFD-fBiCGSTAB algorithm.
	\begin{table}[!htbp]\small
		\centering \caption{Performance of the fractional CN-BCFD schemes  solved by different solvers  with fixed $N=2^{12},\alpha=1.8,\ \kappa=1$}
		\setlength{\tabcolsep}{1.8mm}\label{tab:1d-1}
		\begin{tabular}{l l l l l }
			\toprule
			Algorithm	&M          & Error-u       &CPU             & Itr.  \\
			\midrule
			&$2^7$      &1.1985e-02       &4.63s              &             \\
			&$2^8$      &6.7979e-03       &39s                &             \\
			CN-BCFD-GE&$2^9$      &3.7886e-03   &692s=11m 32s         &             \\
			&$2^{10}$   &2.0818e-03   &5189s=1h 26m 29s       &             \\
			&$2^{11}$   &1.1304e-03   &166821s= 46h 20m 21s     &             \\
			\midrule
			&$2^7$      &  1.1985e-02    & 0.64s             &    2.00     \\
			&$2^8$      &  6.7979e-03   & 2.86s               &    2.65      \\
			CN-BCFD-BiCGSTAB&$2^9$      &  3.7886e-03 & 20.15s        &    5.06   \\
			&$2^{10}$   &  2.0818e-03 & 171s    &   9.67     \\
			&$2^{11}$   &  1.1304e-03 & 4626s=1h 17m 6s &    18.97    \\\midrule
			&$2^7$      & 1.1985e-02     & 1.10s   &    2.00      \\
			&$2^8$      & 6.7979e-03      & 2.56s   &   2.64   \\
			CN-BCFD-fBiCGSTAB&$2^9$      & 3.7886e-03  & 8.71s   &    5.06   \\
			&$2^{10}$   & 2.0818e-03  & 31.83s  &    9.68  \\
			&$2^{11}$   & 1.1304e-03  & 157s=2m 37s    &     18.94   \\      
			\bottomrule
		\end{tabular}
	\end{table}	
	\begin{table}[htbp]\small
		\centering \caption{Performance of the fractional CN-BCFD schemes solved by different solvers  with fixed $N=2^{12}$, $\alpha=1.6$, $\kappa=1.5$}
		\setlength{\tabcolsep}{1.8mm}\label{tab:1d-2}
		\begin{tabular}{l l l l l  }
			\toprule
			Algorithm	&M          & Error-u     &CPU                  & Itr.  \\\midrule
			&$2^7$      &2.9450e-03      & 4.49s    &              \\
			&$2^8$      &1.1803e-03      & 40s             &             \\
			CN-BCFD-GE&$2^9$      &4.7160e-04      &693s=11m 33s        &             \\
			&$2^{10}$   &1.8695e-04      &5180s=1h 26m 20s     &             \\
			&$2^{11}$   & 7.3360e-05     &197474s=54h 51m 14s &             \\
			\midrule
			&$2^7$      &2.9450e-03      &0.46s             &2.00          \\
			&$2^8$      &1.1803e-03      &2.21s             &2.00          \\
			CN-BCFD-BiCGSTAB&$2^9$     &4.7158e-04     &9s      &2.00          \\
			&$2^{10}$   &1.8695e-04      &56s       &3.51         \\
			&$2^{11}$   &7.3333e-05      &2540s=42m 20s   &5.50          \\
			\midrule
			&$2^7$      &2.9450e-03      &3.15s             &2.00          \\
			&$2^8$      &1.1803e-03      &5.98s             &2.00          \\
			CN-BCFD-fBiCGSTAB&$2^9$    &4.7158e-04    &12s        &2.00          \\
			&$2^{10}$   &1.8695e-04      &35s              &3.53         \\
			&$2^{11}$   &7.3329e-05      &93s=1m 33s              &5.31          \\      
			\bottomrule
		\end{tabular}
	\end{table}
	\subsection{Numerical results for two-dimensional SFDEs}
	\begin{example}\label{exam:2d-1}
		In this example, we extend our discussion to the two-dimensional case \eqref{model-2d:e1}. Let $K^{x,L}=K^{x,R}=K^{y,L}=K^{y,R}=1$, and $f(x,y,t)$ is properly chosen such that the true solution of model \eqref{model-2d:e1} is given by
		\[
		u(x,y,t)=e^t x^2(2-x)^2y^2(2-y)^2,\quad x\in[0,2]^2,\ t\in[0,1].
		\]
	\end{example}
	\begin{table}[htbp]\small
		\centering \caption{Performance of the fractional CN-BCFD schemes solved by different solvers for the two-dimension SFDE with different $\alpha$, $\beta$ and fixed $N=2^{10}, \kappa=1$}\label{tab:2d-1}
		\begin{tabular}{c c c c c c c}
			\toprule
			\multirow{2}{*}{$(\alpha,\beta)$} & \multirow{2}{*}{$M_x=M_y$} & \multicolumn{2}{c}{CN-BCFD-GE} & \multicolumn{3}{c}{CN-BCFD-fBiCGSTAB} \\ 
			\cmidrule(r){3-4}     \cmidrule(r){5-7}     
			&         & Error-u    & CPU  & Error-u     & CPU       & Itr. \\	
			\midrule
			&$2^5$    & 2.7164e-02 & 34m 15s & 2.2680e-02 & 1m 12s  & 3.00 \\
			$(1.5,1.5)$ &$2^6$    & 1.1978e-02 & $\approx $7d & 1.1884e-02 & 4m 30s  & 4.00  \\
			&$2^7$    & ---  &  ---& 6.0684e-03 & 12m 35s & 5.00  \\ 
			&$2^8$    &---  & --- & 3.0603e-03 & 46m 59s & 8.01  \\
			\midrule
			&$2^5$    & 2.1504e-02 & 30m 11s & 2.0971e-02 & 1m 20s & 3.00 \\
			$(1.4,1.7)$ &$2^6$    & --- & --- & 1.1402e-02 & 4m 18s &  5.00 \\
			&$2^7$    & --- & --- & 6.0605e-03 & 18m 12s &  8.00 \\ 
			&$2^8$    & --- & --- & 3.1988e-03 & 1h 22m 49s &  14.00 \\
			\midrule
			&$2^5$    & 1.6065e-02 & 21m 25s & 1.6006e-02 & 1m 40s & 4.00 \\
			$(1.6,1.9)$ &$2^6$    & --- & --- & 8.4224e-03 & 6m 1s &  7.00 \\
			&$2^7$    & --- & --- & 4.3937e-03 & 46m 13s & 14.00  \\ 
			&$2^8$    & --- & --- & 2.2723e-03 & 2h 37m 14s &  27.14 \\
			\bottomrule
		\end{tabular}
	\end{table}
	
	In this simulation, we fix $N=2^{10}$ and test different combinations of $(\alpha, \beta)$ and $\kappa$ for the comparisons of the developed BCFD schemes for the two-dimensional model \eqref{model-2d:e1}. 
	In Tables \ref{tab:2d-1}--\ref{tab:2d-2}, we present the errors and CPU times for the CN-BCFD-GE algorithm and the fast version CN-BCFD-fBiCGSTAB algorithm. Meanwhile, the average numbers of iterations per time level for the latter one are also list. As the one-dimensional case, it is seen that both algorithms generate the same accurate numerical solutions whether on uniform or nonuniform grids. Besides, it also shows that the usage of graded grids evidently improves the computational accuracy, see comparisons of the results in Tables \ref{tab:2d-1}--\ref{tab:2d-2}. In addition, the fast version CN-BCFD-fBiCGSTAB algorithm shows an efficient and effective capability for solving the two-dimensional model. For example, with $(\alpha,\beta)=(1.5,1.5)$, $\kappa=1$, and $M_x=M_y=2^6$, the CN-BCFD-GE algorithm spends around 7 days to reach the error of magnitude $10^{-2}$ , while the CN-BCFD-fBiCGSTAB algorithm consumes only no more than 5 minutes for the nearly same error. However, for larger-scale modeling, for example, $(\alpha,\beta)=(1.5,1.5)$, $M_x=M_y=2^8$ and $\kappa=1.5$, the CN-BCFD-GE algorithm can hardly run on the computer, but the CN-BCFD-fBiCGSTAB algorithm can still achieve the error of  magnitude $10^{-4}$ in only three and a half  hours. It is believed that the comparisons shall be more obvious for even larger $M_x$ and $M_y$. Finally, we would like to point out that with the increasing number of grids, the average number of iterations becomes larger and larger especially on the nonuniform grids, see Table  \ref{tab:2d-2}. This in turn increases the overall computational cost. Therefore, it is of great importance to develop an efficient preconditioning technique for the resulting linear algebraic system \cite{DJWWZ21,FSW20,FLW19,FW17,JW15,LFWC19,PKNS14,WD13,XLS22}. 
	To sum up, numerical results for the two-dimensional model further demonstrate the strong superiority of the CN-BCFD-fBiCGSTAB algorithm in large-scale modeling and simulation. 
	\begin{table}[!htbp]\small
		\centering \caption{Performance of the fractional CN-BCFD schemes solved by different solvers for the two-dimension SFDE with different $\alpha$, $\beta$ and fixed $N=2^{10}, \kappa=1.5$}\label{tab:2d-2}
		\begin{tabular}{c c c c c c c}
			\toprule
			\multirow{2}{*}{$(\alpha,\beta)$} & \multirow{2}{*}{$M_x=M_y$} & \multicolumn{2}{c}{CN-BCFD-GE} & \multicolumn{3}{c}{CN-BCFD-fBiCGSTAB} \\ 
			\cmidrule(r){3-4}     \cmidrule(r){5-7}     
			&         & Error-u    & CPU  & Error-u     & CPU       & Itr.  \\	
			\midrule
			&$2^5$    & 1.8287e-02 & 46m 49s & 1.8302e-02 & 2m 3s  & 4.66 \\
			$(1.5,1.5)$ &$2^6$    & --- & --- & 4.9508e-03 & 9m 8s  & 7.06  \\
			&$2^7$    & ---  &  ---& 2.1587e-03 & 37m 43s & 13.90  \\ 
			&$2^8$    &---  & --- & 9.9896e-04 & 3h 33m 38s & 31.20  \\
			\midrule
			&$2^5$    & 1.7733e-02 & 21m26s & 1.7752e-02 & 2m 20s & 5.00 \\
			$(1.4,1.7)$ &$2^6$    & --- & --- & 4.7589e-03 & 12m 33s &  10.00 \\
			&$2^7$    & --- & --- & 1.6591e-03 & 1h 4m 16s &  23.41 \\ 
			&$2^8$    & --- & --- & 7.9748e-04 & 7h 20m 47s &  58.49 \\
			\midrule
			&$2^5$    & 1.7347e-02 & 21m 24s & 1.7359e-02 & 3m 9s & 7.01 \\
			$(1.6,1.9)$ &$2^6$    & --- & --- & 4.4954e-03 & 19m 12s &  17.67 \\
			&$2^7$    & --- & --- & 1.1314e-03 & 2h 18m 7s & 47.62  \\ 
			&$2^8$    & --- & --- & 2.7873e-04 &16h 31m 20s &  139.67 \\
			\bottomrule
		\end{tabular}
	\end{table}
	
	\section{Conclusions}
	This paper seems to be the first attempt to discuss fractional BCFD method and its efficient implementation for both one and two-dimensional SFDEs with fractional Neumann boundary condition on general nonuniform grids. Totally different from the case of uniform spatial partitions, where Toeplitz-like matrix structure can be discovered for the widely used numerical methods such as finite difference and finite volume methods, and that would lead to significantly reduced memory requirement and computational complexity.  The developed fractional BCFD method is well suitable for the nonuniform grids, and based upon the SOE approximation technique, a fast version fractional BCFD algorithm is further constructed, in which fast matrix-vector multiplications of the resulting coefficient matrices with any vector in only $\mathcal{O}(MN_{exp})$ operations are developed. Numerical experiments show that the usage of locally refined grids greatly improves the accuracy of the developed method, and most importantly, the fast version fractional CN-BCFD algorithm has greatly reduced the CPU times without any accuracy being lost. 
	
The method shows strong potential in large-scale modeling and simulation on nonuniform spatial grids. We will continue studying this interesting topic at the next step, for example, combinations with the dimensional splitting method maybe  is a good idea for simulating high-dimensional SFDEs \cite{WD13_3D}. 
However, due to the absence of regularity information on the solution, conducting an error analysis for the proposed scheme remains a challenging task. Moreover, we still have no idea about the unique solvability of the solution to the fractional BCFD method on general nonuniform grids, which also presents a significant challenge for future work. 
	
\begin{acknowledgements}
The authors would like to express their most sincere thanks to the referees for their very helpful comments and suggestions, which greatly improved the quality of this paper. This work was supported in part by the National Natural Science Foundation of China (No. 11971482), by the Fundamental Research Funds for the Central Universities (No. 202264006) and by the OUC Scientific Research Program for Young Talented Professionals. 

\end{acknowledgements}

\vspace{-0.6cm}
	\section*{Declarations}
	{\small \textbf{Conflict of interest}	The authors declare that they have no conflict of interest.}

	
	
	\bibliographystyle{spmpsci}    
	\bibliography{myref}

	\bigskip  
	
	\small 
	\noindent
	{\bf Publisher's Note}
	Springer Nature remains neutral with regard to jurisdictional claims in published maps and institutional affiliations.
\end{document}